\documentclass[11pt,final]{amsart}
\usepackage{amsmath, amsthm, amsfonts, ifpdf}
\usepackage{color}
\usepackage{graphicx}
\usepackage{multicol}
\usepackage{epstopdf}

\usepackage{tikz}
\usetikzlibrary{arrows,shapes,trees,backgrounds}
\usetikzlibrary{intersections}

\theoremstyle{plain}
\newtheorem{theorem}{Theorem}[section]
\newtheorem*{theorem*}{Theorem}

\newtheorem{pro}[theorem]{Proposition}
\newtheorem{Def}[theorem]{Definition}
\newtheorem{lem}[theorem]{Lemma}
\newenvironment{customthm}[1]
{\renewcommand\theinnercustomthm{#1}\innercustomthm}
{\endinnercustomthm}
\newtheorem{innercustomthm}{Theorem}

\newtheorem{cor}[theorem]{Corollary}

\theoremstyle{definition}
\newtheorem*{Def*}{Definition}

\newtheorem{Rem}[theorem]{Remark}

\numberwithin{equation}{section}

\newcommand{\bpo}{\begin{pro}}
\newcommand{\epo}{\end{pro}}
\newcommand{\be}{\begin{equation}}
\newcommand{\ene}{\end{equation}}
\newcommand{\br}{\begin{Rem}}
\newcommand{\er}{\end{Rem}}
\newcommand{\bl}{\begin{lem}}
\newcommand{\el}{\end{lem}}
\newcommand{\bd}{\begin{Def}}
\newcommand{\ed}{\end{Def}}
\newcommand{\ben}{\begin{enumerate}}
\newcommand{\een}{\end{enumerate}}
\newcommand{\bp}{\begin{proof}}
\newcommand{\ep}{\end{proof}}
\newcommand{\beq}{\begin{equation*}}
\newcommand{\eeq}{\end{equation*}}
\newcommand{\bear}{\begin{eqnarray*}}
\newcommand{\eear}{\end{eqnarray*}}
\newcommand{\bt}{\begin{theorem}}
\newcommand{\et}{\end{theorem}}
\newcommand{\bst}{\begin{split}}
\newcommand{\est}{\end{split}}

\newcommand{\bal}{\begin{aligned}}
\newcommand{\eal}{\end{aligned}}

\renewcommand{\P}{\partial}
\newcommand{\F}[2]{\frac{#1}{#2}}
\newcommand{\la}{\langle}
\newcommand{\ra}{\rangle}

\newcommand{\R}{\mathbb{R}}

\newcommand{\bnb}{\bar{\nabla}}
\newcommand{\nb}{\nabla}

\newcommand{\Sc}{\varepsilon}

\newcommand{\vp}{\varphi}

\newcommand{\tw}{\tilde{\omega}}
\newcommand{\PLH}{{\mkern-1mu\times\mkern-1mu}}

\newcommand{\TIMES}{\times}

\newcommand{\gu}{graph$(u)$}

\def\XXint#1#2#3{{\setbox0=\hbox{$#1{#2#3}{\int}$}
    \vcenter{\hbox{$#2#3$}}\kern-.5\wd0}}

\makeatletter
\def\@citestyle{\m@th\upshape\mdseries}
\def\citeform#1{{\bfseries#1}}
\def\@cite#1#2{{%
  \@citestyle[\citeform{#1}\if@tempswa, #2\fi]}}
\@ifundefined{cite }{%
  \expandafter\let\csname cite \endcsname\cite
  \edef\cite{\@nx\protect\@xp\@nx\csname cite \endcsname}%
}{}
\makeatother
	
\begin{document}
		\title[Nonparametric mean curvature flow with angle conditions]{Nonparametric mean curvature type flows of graphs with contact angle conditions}
		\author{Hengyu Zhou}
		\date{\today}
		\address{Department of Mathematics, Sun Yat-sen University, No. 135, Xingang Xi Road, Guangzhou, 510275, P. R. China}
		\email{hyuzhou84@yahoo.com}
		\date{\today}
	\begin{abstract}
		In this paper we study nonparametric mean curvature type flows in $M\PLH\R$ which are represented as graphs $(x,u(x,t))$ over a domain in a Riemannian manifold $M$ with prescribed contact angle. The speed of $u$ is
		the mean curvature speed minus an admissible function $\psi(x,u,Du)$. \\
		\indent Long time existence and uniformly convergence are established if $\psi(x,u, Du)\equiv 0$ with
		vertical contact angle and $\psi(x,u,Du)=h(x,u)\omega$ with $h_u(x,u)\geq h_0>0$ and $\omega=\sqrt{1+|Du|^2}$. Their applications include mean curvature type equations with prescribed contact angle boundary condition and the asymptotic behavior of nonparametric mean curvature flows of graphs over a convex domain in $M^2$ which is a surface with nonnegative Ricci curvature.
	\end{abstract}
\maketitle
\section{Introduction}\label{sec:intro}
The main topic of this paper is various nonparametric mean curvature type flows of graphs with prescribed contact angle in product manifolds. These flows can be viewed as a generalization of nonparametric mean curvature flows. They are deeply rooted in mean curvature type equations with prescribed contact angle boundary condition to which certain uniformly convergence of these flows may yield the solution.\\
\indent Throughout this paper $M$ denotes a closed Riemannian manifold with a metric $\sigma$ and $\Omega$ is a smooth bounded domain of $M$. A nonparametric mean curvature type flow here is a family of graphs $(x,u(x,t))$ for $x\in \Omega$ in $M\PLH \R$ where $u(x,t):\Omega \times [0,\R)\rightarrow \R$ is the solution to the following evolution equation
\be\label{eq:first}
\begin{aligned}
	&u_t =g^{ij}(Du)u_{ij}-\psi(x,u,Du)\quad g^{ij}(Du)u_{ij}=div(\F{Du}{\omega})\omega \\
	&u_{\gamma} =\phi(x)\omega  \quad \text{on $\P\Omega\times \ [0,\infty)$} \quad u(x,0)=u_0(x) \quad  x\in \bar{\Omega};
\end{aligned}
\ene
where $Du$ is the gradient of $u(x,t)$, $\omega$ denotes $\sqrt{1+|Du|^2}$ and $div$ is the divergence of $M$.
Here $\gamma$ denotes the interior normal to $\P\Omega$ and $\phi(x) \in C^{\infty}(\bar{ \Omega})$ satisfies the property $|\phi|\leq \phi_0<1$. The lower index $i,j,\gamma$ indicate covariant derivatives with respect to the metric $\sigma$.
We call $g^{ij}(Du)u_{ij}$ as the \emph{mean curvature speed} since
$$
g^{ij}(Du)u_{ij}=div(\F{Du}{\omega})\omega=H\omega
$$
where $H$ is the mean curvature of the graph of $u(x)$ in $M\PLH\R$. \\
\indent We are also interested in the solution of a corresponding elliptic version of problem \eqref{eq:first} as follows.
\be\label{eq:second}
\begin{aligned}
	g^{ij}(Du)u_{ij}&=\psi(x,u,Du) &\quad &\text{in $\Omega$} \\
	u_{\gamma} &=\phi(x)\omega &\quad  \quad &x\in \P\Omega;
\end{aligned}
\ene
\indent Problems \eqref{eq:first} and \eqref{eq:second} have been studied extensively in various particular settings. In the case of $\psi\equiv 0$ problem \eqref{eq:first} describes a \emph{nonparametric mean curvature flow}. It is a family of graphs over a domain with prescribed contact angle flowed by mean curvature vector in product manifolds. General long time existence of the solution and uniformly convergence were investigated by \cite{Hui89A}, \cite{GB96} and recently \cite{Xu16} in Euclidean spaces. The graph of $u_\infty(x)+Ct$ in $M\PLH \R$ is called a translating solution to mean curvature flow with prescribed contact angle if $u_\infty(x)$ is a solution to problem \eqref{eq:second} with $\psi(x,u,Du)\equiv C$ (see equation \eqref{eq:dc}). The translating surface in $\R^3$ was investigated by \cite{AW94}. When $\psi(x,u, Du)$ is $h(x,u)\omega$ with $h_u(x,u)\geq h_0>0$, problem \eqref{eq:second} are well known as \emph{Capillary problems with positive gravity} (see \cite{CF74}). Various existence results of their solutions have obtained by \cite{Ger76}, \cite{Spr75}, \cite{Ura80},\cite{Ko88},\cite{Lie88} in Euclidean spaces. \\
\indent The main goal of this paper is to  develop a unified framework about mean curvature type flows with prescribed contact angle and apply this to reexamine Capillary problems with positive gravity. Our results shall include all previous results in this direction (mainly in \cite{Hui89A,GB96,AW94}). The main obstacle is the curvature of $M$. We remind the reader that even in the case of mean curvature flows ($\psi\equiv 0$) there is an essential difference between those in Eucldean spaces and general Riemannian manifolds (see \cite{Hui86}).  \\
\indent We need  a concept regarding $\psi(x,u, Du)$ such that problems \eqref{eq:first} and \eqref{eq:second} make sense. It is said that $\psi(x,u,Du)$ is an \emph{admissible} function if there is a constant $C\geq 0$ with the following property:
\be\label{condition:c1a}
\tag{c1}\left\{\begin{split}
	&\psi_{u}(x, u, Du)\geq -C,\quad    \sigma^{kl}\psi_k(x, u, Du)\psi_l(x, u, Du)\leq C\omega^2 \\
	& \psi(x,u,Du)-\psi_{u_k}(x,u, Du)u_k\geq -C, \quad   \sigma_{kl}\psi_{u_k}\psi_{u_l}\leq C
\end{split}\right.
\ene
where $\sigma=\sigma_{ij}dx_i dx_j$ on $M$ and $(\sigma^{ij})$ is its inverse of $(\sigma_{ij})$. More detail is given in page \pageref{condition:c1}. Two examples of admissible functions are $h(x,u)\omega$ with $h_u\geq 0$ and $Bu$ for a constant $B$. See Lemma \ref{lm:condition}.\\
\indent Next we explain our idea to establish the gradient estimate of the solution in problem \eqref{eq:first}. Let $u(x,t)$ be the solution to problem \eqref{eq:first} for $\psi(x,u, Du)$ with an admissible constant $C$. We construct an auxiliary function
$$
\eta = e^{Ku(x,t)}(Nd(x)+1-\phi(x) \la \vec{v}, Dd\ra )
$$
where $d(x)$ is a smooth function and equal to $d(x,\P\Omega)$ provided $x$ is sufficiently close to $\P\Omega$, $\vec{v}$ is the downward normal to the graph of $u(x,t)$. Two positive constants $N$ and $K$ are determined as follows. Given any $T>0$ suppose that the maximum of $\eta\omega$ on $\bar{\Omega}\times [0,T]$ is achieved at $(x_0,t_0)$, then
\begin{itemize}
	\item [(i)]  we can choose a $N$, sufficiently large and independent of $K$, such that $\omega(x_0,t_0)\leq K$ whenever $x_0\in \P\Omega$(see Lemma \ref{lemma:N});
	\item [(ii)] after fixing this $N$, we can choose sufficiently large $K$ such that $\omega(x_0,t_0)\leq C_0$ whenever $x_0\in \Omega$. Here $C_0$ relies on $d(x), \phi(x)$, the admissible constant $C$ and the Ricci curvature on $\bar{\Omega}$ (see Lemma \ref{lm:bounda}).
\end{itemize}
Notice that the choice of $N$ and $K$ are also independent of $T$ and the position of $(x_0, t_0)$. This method is inspired by \cite{GB96} which dealt with problem \eqref{eq:first} in Euclidean spaces. A remarkable point is that the Ricci curvature of $M$ is absorbed when $K$ is chosen. A direct derivation in Theorem \ref{thm:boundge} yields a crucial estimate in problem \eqref{eq:first}. That is $\omega\leq C_0 e^{2K U_T}$ where $C_0$ is a constant only depending on $C,d(x),\phi(x)$ and $U_T$ is the maximum of $|u(x,t)|$ on $\bar{\Omega}\times [0, T]$.  \\
\indent We are ready to state main results in this paper. First we see that the long time existence for solutions in problem \eqref{eq:first} is a very general fact.
\bt\label{thm:GLT} Suppose $\psi(x,u,Du)$ is a smooth function of $x,u$ and $Du$ satisfying 
\begin{enumerate}
	\item [(a)]$\psi_u(x,u,Du)\geq c_0$ where $c_0$ does not rely on $u$;
	\item [(b)]for any given $K>0$ there is a positive constant $C=C(K)$ such that $\psi(x,u,Du)$ is admissible with respect to $C$ if $|u|\leq K$.
\end{enumerate}  Then problem \eqref{eq:first} with $\psi(x,u,Du)$ has a solution $u\in C^{\infty}(\bar{\Omega}\times [0,\infty))$.
\et
It extends Theorem 2.4 in \cite{GB96}.\\
\indent The following result is very natural as a generalization of the main result in \cite{Hui89A}.
\bt\label{thm:MT1}  Problem \eqref{eq:first} with $\psi\equiv \phi\equiv 0$ has a smooth solution $u(x,t)$ on $\bar{\Omega}\times[0, \infty)$. Moreover $u(x,t)$ converges uniformly to a constant as $t$ goes to infinity.
\et
The determination of this constant may have independent interests. 
Notice that the graph of $u(x,t)$ is a (nonparametric) mean curvature flow in the product manifold $M\PLH \R$. But the product structure of $M\PLH \R$ can not be weakened as a warped product structure. In appendix A we will construct a graph in a warped product manifold. Its mean curvature flow will break the graphical property and form a singularity in finite time. In this sense Theorem \ref{thm:MT1} is the best result we can expect.\\
\indent The following result is inspired by the first example in (\cite{GB96}, Section 3). We discover a general connection between Capillary problems with positive gravity and nonparametric mean curvature flows.
\bt\label{thm:MT2} Let $\psi(x,u, Du)$ be $h(x,u)\omega$ with $h_u\geq h_0>0$. Problem \eqref{eq:first} has a smooth solution $u(x,t)$ on $\bar{\Omega}\times [0,\infty)$. Moreover $u(x,t)$ converges uniformly to $u_\infty(x)$ as the solution to problem \eqref{eq:second} as $t$ goes to infinity.
\et
A straightforward consequence of Theorem \ref{thm:MT2} is that problem \eqref{eq:second} has a unique smooth solution when $\psi(x,u,Du)$ is $h(x,u)\omega$ with $h_u\geq h_0>0$ (see Theorem \ref{thm:MT22}). \\
\indent A common feature of Theorem \ref{thm:MT1} and Theorem \ref{thm:MT2} is that the uniform bound of $u(x,t)$ indicates that of $\omega$. However this is not the only model to characterize the evolution in problem \eqref{eq:first}. The following result, a generalization of the work in \cite{AW94}, describes a different one.
\bt\label{thm:MT3} Suppose $M^2$ is a Riemannian surface with nonnegative Ricci curvature. Let $\Omega$ be a bounded convex domain with
$$k_0 \geq k\geq |\F{\phi_T}{\sqrt{1-\phi^2}}|+ \delta_0$$
where $k$ is the inward curvature of $\P\Omega$ and $T$ is the tangent vector of $\P\Omega$, $k_0,\delta_0$ are positive constants.
Then
\begin{itemize}
	\item [(i)]problem \eqref{eq:first} with $\psi(x,u,Du)\equiv 0$ has a smooth solution $u(x,t)$ on $\bar{\Omega}\times[0,\infty)$ with  $\omega\leq C_1$ where $C_1$ is a constant depending on $ k_0, \delta_0, \phi_0$ and $u_0(x)$.
	\item [(ii)]
	moreover $u(x,t)$ converges uniformly to $u_{\infty}(x)+Ct$ as $t\rightarrow \infty$ where $u_{\infty}(x)$ is the solution to problem \eqref{eq:second} with $\psi(x,u, Du)\equiv C$. Here $C$ is given by
	$$
	C=-\F{\int_{\P\Omega}\phi(x)dx}{\int_\Omega (1+|Du_{\infty}|^2)^{-\F{1}{2}}dx}
	$$
\end{itemize}
\et
We remark that the existence of $u_{\infty}(x)$ above is highly nontrivial. It follows from four conditions: the solvability of Capillary problems of positive gravity, the nonnegative Ricci curvature, the convex assumption about $\P\Omega$ and the fact that $M^2$ is a surface. See the proof of Lemma \ref{lm:ts}.\\
\indent
The paper is organized as follows. In section 2 we discuss preliminary facts for later reference. Section 3 is devoted to several estimates along the flow in problem \eqref{eq:first}. In section 4 we establish Theorem \ref{thm:GLT}, Theorem \ref{thm:MT1} and Theorem \ref{thm:MT2} with the estimates in section 3. A key point is to control $u(x,t)$ under various settings. We give a sufficient condition about the uniformly convergence of the flow in problem \eqref{eq:first} in Lemma \ref{lm:cn}. In section 5 we treat asymptotic behaviors of nonparametric mean curvature flows over a convex domain in Riemannian surfaces and obtain Theorem \ref{thm:MT3}. In appendix A we construct a graph in a warped product manifold. Its mean curvature flow only exists smoothly in finite time. This part is independent of other parts in this paper. In appendix B we give an area formula for graphs in warped product manifolds. In appendix C some technique lemmas about parabolic equations are presented. In appendix D we show Lemma \ref{lm:cn}.

\section{The Geometry of Graphs}
\indent In this section we present some preliminary facts of graphs in $M\PLH \R$. Then we discuss the admissible condition of $\psi(x,u, Du)$.\\
\indent Suppose $u(x)$ is a smooth function on $\bar{\Omega}$. Recall that $\sigma$ is the Riemannian metric on $M$. A product manifold $M\PLH\R$ is the set $\{(x,r):x\in M, r\in \R\}$ equipped with the metric $\sigma+dr^2$. The graph of $u(x)$ in $M\PLH \R$ is denoted by {\gu}.\\
\indent Vectors on $M$ will be denoted by $\{X^i\}$ , covectors by $\{Y_i\}$ and mixed tensors by $T=\{T^{i}_{jk}\}$. We always sum over repeated indices from $1$ to $n$ and use brackets for inner products on $M$. Consider the inner product of mixed tensors given as follows.
$$
\la T, S\ra =\la T_{jk}^i, S_{jk}^i\ra =\sigma_{is}\sigma^{jr}\sigma^{kt}T^{i}_{jk}S^s_{rt},\quad |T|^2=\la T^i_{jk}, T^i_{jk}\ra
$$
A Cauchy inequality is stated as
\be\label{Cauchy:Inequality}
| \la T, S\ra|\leq |T||S|
\ene
Choose a local coordinate $\{x_i\}$ on $\Omega$ and $r$ on $\R$ respectively. We write $\F{\P}{\P x_i}, \F{\P}{\P r}$ as $\P_i, \P_r$ for short. Let $u_i, u_{ij}$ be the covariant derivatives of $u$ with respect to $\P_i,\P_j$. We collect the following notation.
\begin{gather}
\sigma_{ij}= \la \P_i, \P_j\ra, \quad (\sigma^{ij})=(\sigma_{ij})^{-1},\quad u^i=\sigma^{ik}u_k\\
Du=u^i\P_i,\quad \omega=\sqrt{1+|Du|^2}\quad \vec{v}=\F{1}{\omega}(Du -\P_r)\label{def:vector}
\end{gather}
where $Du$ is the gradient of $u$ and $\vec{v}$ is the downward normal unit vector to {\gu}. Thus
$\{X_i:\P_i+u_i\P_r\}_{i=1}^n
$ is a local frame of {\gu}.\\
\indent
The first fundamental form of {\gu} is
$$g_{ij}=\la X_i, X_j\ra=\sigma_{ij}+u_iu_j $$
with the inverse
\be
g^{ij}(Du)=\sigma^{ij}-\F{u^iu^j}{\omega^2}
\ene
A direct computation gives $g^{ij}(Du)u_{ij}=div(\F{Du}{\omega})\omega$.  The following identities and inequalities will be used in the proof of Lemma \ref{lm:bounda}.
\bl The partial derivatives of $g^{ij}$ fulfill that
\be \label{parder:g}
\F{\P g^{ij} }{\P u_k}=-\F{1}{\omega^2}(g^{ik}u^j+g^{jk}u^i)\quad\quad \F{\P g^{ij}}{\P u_k} u_{ij}=-\F{2}{\omega}g^{ik}\omega_i
\ene
where $\omega_i=\F{u^ku_{ki}}{\omega}$. Consider three mixed tensors $T=\{T_{ij}\}$, $S=\{S_i\}$ and $U=\{U_i\}$. Then it follows that
\be \label{estimate:G}
|g^{ij}T_{ij}|\leq (n+1)|T|, \quad |g^{ij} S_i U_j|\leq 2|S||T|
\ene
where $n$ is the dimension of $M$.
\el
\bp Equation \eqref{parder:g} follows from a direct computation.  Expanding $g^{ij}T_{ij}$ and applying the Cauchy inequality we observe that
\begin{align*}
|g^{ij}T_{ij}|&= |\sigma^{ij}T_{ij}-\F{1}{\omega^2}u^iu^jT_{ij}|\\
&\leq n|T|+|\F{1}{\omega^2}\la u_iu_j, T_{ij}\ra| \leq (n+1)|T|
\end{align*}
As for $g^{ij}S_i U_j$, a similar derivation yields that
\begin{align*}
|g^{ij}U_iS_j|&\leq |\sigma^{ij}U_iS_j-\F{1}{\omega^2}u^iu^jU_iS_j|\\
&\leq |\la U_i, S_i\ra| +|\F{1}{\omega^2}\la u_i u_j, U_i S_j\ra |\leq 2|U||S|
\end{align*}
We obtain the desired conclusion.
\ep
The second fundamental form of {\gu} is
\be
h_{ij}=\la \bnb_{X_i}\vec{v}, X_j\ra =\F{u_{ij}}{\omega};
\ene
where $\bnb$ is the covariant derivative of $M\PLH \R$. The norm of the second fundamental form is
\be
|A|^2=g^{kl}g^{ij}h_{ki}h_{lj} =\F{g^{il}g^{kj}u_{kl}u_{ij}}{\omega^2};
\ene
We write $-H\vec{v}$ for the mean curvature vector of {\gu} where $H=g^{ij}h_{ij}$. Thus we see that
\be
g^{ij}(Du)u_{ij}=H\omega=div(\F{Du}{\omega})\omega
\ene
This is the reason that we call the solution to problem \eqref{eq:first} as a mean curvature type flow. Moreover
\bl Let $u(x,t)$ be a solution to problem \eqref{eq:first} with $\psi(x,u, Du)\equiv 0$. Then the graph $F(x,t)=(x, u(x,t))$ fulfills that
\be (\F{d F}{dt})^{\bot}=-H\vec{v}
\ene
where $\bot$ is the projection into the normal bundle of $F(.,t)$ in $M\PLH\R$.
\el
To study problem \eqref{eq:first} it is useful to propose a concept about $\psi(x,u,Du)$ for later reference. This concept is inspired by the assumptions of Lemma 2.3 in \cite{GB96}.\\
\indent
Given two smooth functions $u(x)$ and $\psi(x,u,Du)$ we define a smooth vector and a smooth covector along $\bar{\Omega}$ as follows
$$
X_\psi=\psi_{u_k}\P_k,\quad Y^\psi = \psi_l dx^l
$$
where  $\psi_k:=\F{\P\psi(x, u, Du)}{\P x_k}$ and $\psi_{u_k}:=\F{\P \psi}{\P u_k}(x,u,Du)$. These definitions are well-defined. A straightforward verification yields that $X_\psi$ and $Y^\psi$ are independent of local coordinates.  \\
\indent We say that $\psi(x,u,Du)$ is \emph{admissible} with respect to a constant $C\geq 0$ if it fulfills that
\be\label{condition:c1}
\tag{c1}\left\{\begin{split}
	&\psi_{u}(x, u, Du)\geq -C,\quad    \sigma^{kl}\psi_k(x, u, Du)\psi_l(x, u, Du)\leq C\omega^2 \\
	& \psi(x,u,Du)-\psi_{u_k}(x,u, Du)u_k\geq -C, \quad   \sigma_{kl}\psi_{u_k}\psi_{u_l}\leq C
\end{split}\right.
\ene
which is equivalent to
\be
\left\{\begin{split}
	&\psi_{u}(x, u, Du)\geq -C,\quad    \la Y^\psi, Y^\psi\ra \leq C\omega^2 \\
	& \psi(x,u,Du)-du(X_\psi)\geq -C, \quad   \la X_\psi,X_\psi\ra\leq C
\end{split}\right.
\ene
Hence our definition is well-defined.
The conditions in assumptions \eqref{condition:c1} are very general. Some admissible functions are listed in the following result.
\bl \label{lm:condition} Suppose $h(x,u)$ is smooth. If a smooth function $u$ satisfies $|u|\leq K$ and one of the followings holds:
\begin{itemize}
	\item [(i)] $\psi(x,u,Du)$ is $h(x,u)\omega$ with $h_u(x,u)\geq 0$;
	\item  [(ii)]$\psi(x,u,Du)$ is $h(x,u)$ with $h_{u}(x,u)\geq h_0$ for some constant $h_0$,
\end{itemize}
where $h_u=\F{\P h(x,u)}{\P u}$, then $\psi(x,u,Du)$ is admissible with respect to some constant $C$ only depending on $K$ or $K$ and $h_0$.
\el
\bp The proof is a straightforward computation.
\ep
Let $\nb$ be the covariant derivative of $M$. The Riemann curvature tensor of $M$ is
\be
R(X,Y)Z=\nb_{Y}\nb_{X}Z-\nb_{X}\nb_{Y}Z+\nb_{[X,Y]}Z
\ene
Then we write $ R(X,Y,Z,W)$ for $\la R(X,Y)Z, W\ra$. Now we introduce a commuting formula for covariant derivatives.
\bl\label{lm:inchange} Let $\vp_i dx^i$ be a covector on $M$. Then
$$\vp_{ijk}=\vp_{ikj}+R_{kjip}\vp^p $$
where $R_{kijp}=R(\P_k,\P_j,\P_i,\P_p)$ and $\vp^p=\sigma^{pk}\vp_k$.
\el
The proof follows from a straightforward verification. We skip it here.\\
\indent Now we present a lemma regarding the local frame of $\P\Omega$. Recall that $\gamma$ is the interior normal to $\P\Omega$. Notice that $d(x,\P\Omega)$ is a smooth function only for $x$ sufficiently close to $\P\Omega$. In the remainder of this paper, we assume $d(x)$\label{def:dx} is a nonnegative smooth function such that $d(x)=d(x,\P\Omega)$ for $x$ close to $\P\Omega$ and $\la Dd, Dd\ra\leq 1$. The existence of $d(x)$ follows from the geodesic flow of $\P\Omega$ with the initial speed $\gamma$. The following result is easily verified.
\bl \label{eq:local coordinate} Given any point on $\P\Omega$ and one of its small neighborhoods, we can construct a local coordinate $\{x_i\}_{i=1}^{n}$ such that $\{\P_1,\cdots,\P_{n-1}\}$ restricted on $\P\Omega$ is a local frame on $\P\Omega$ and
$$
\la \P_i, \P_n\ra =0,\quad i=1,\cdots, n-1 \quad \la \P_n, \P_n\ra =1
$$
in this small neighborhood of $\P\Omega$. Here $\P_n=Dd(x)$ and $\P_n=\gamma$ on $\P\Omega$.
\el
Together with $d(x)$, $\phi(x)$ and two given positive constants $N,K$ we can define
\be\label{def:eta}
\eta=e^{Ku}(Nd(x)+1-\phi(x) \la \vec{v}, Dd\ra)
\ene
where $u$ is any smooth function on $\bar{\Omega}$ and $\vec{v}$ is the downward normal to {\gu}.
By equation \eqref{def:vector} we can write $\eta$ as

\be\label{eq:def:eta}
\eta=e^{Ku}(Nd(x)+1-\phi(x) v^kd_k)
\ene
where $v^k=\F{u^k}{\omega}$. \\
\indent
One advantage of $\eta$ is that we can control $\omega$ when $\eta\omega$ achieves its maximum on the boundary for appropriate $N$. The following result is inspired by (\cite{GB96}, Lemma 2.2).
\bl\label{lemma:N} Let $u(x)$ be a smooth function on $\bar{\Omega}$ with $u_\gamma=\phi(x)\omega$ on $\P\Omega$ and $\eta$ be given by \eqref{eq:def:eta}. Provided $N$ is sufficiently large, independent of $K$, then
if the maximum of $\eta\omega$ on $\bar{\Omega}$ occurs at $x_0\in\P\Omega$, then
$$
\omega(x_0)\leq K
$$
Here $N$ only depends on $\phi(x)$ and $d(x)$.
\el
\br\label{cordinate} For a fixed point $x_0\in \P\Omega$ we can choose a local coordinate $\{x_i\}_{i=1}^n$ on $M$ satisfying at $x_0$
\begin{itemize}
	\item [(i)]$\{\P_i\}$ is orthonormal and unit, i.e., $\la \P_i,\P_j\ra(x_0)=\delta_{ij}$ and $\nb_{\P_i}\P_j(x_0)=0$ where $\nb$ is the corvariant derivative of $M$,
	\item [(ii)]$\P_n=\gamma$ and
	\be\label{u:assumption}
	u_1\geq 0\quad u_k=0\quad \text{for}\quad 2\leq k\leq n-1;
	\ene
\end{itemize}
\er
\bp Suppose $\eta\omega$ achieves its maximum on $\bar{\Omega}$ at $x_0\in \P\Omega$. Then we take the local coordinate in Remark \ref{cordinate}. From now on all expressions are evaluated at $x_0$.\\
\indent Since $d(x)$ is equal to $d(x,\P\Omega)$ in a neighborhood of $\P\Omega$ then at $x_0$ we have
\be\label{der:d}
d_{\alpha}=0,\quad \text{for}\quad  1\leq \alpha\leq n-1, \quad d_n=1,\quad d_{kn}=0\quad\text{for any $k$}
\ene
The fact $d_{nn}=0$ follows from $\la Dd, Dd\ra=1$ near $x_0$ and $Dd=\gamma=\P_n$ at $x_0$. \\
\indent  Since $\eta\omega$ achieves the maximum on $\bar{\Omega}$ at $x_0$, then
\begin{gather}
\omega_1\eta +\omega\eta_1=0 \label{omega1}\\
0\geq (\omega\eta)_n=e^{Ku}(Ku_n(1-\phi^2)\omega+N\omega+\omega_n-\phi u_{nn}-\phi_nu_n)\label{comomega}
\end{gather}
Here we apply $\eta=e^{Ku}(1-\phi^2)$ on $\P\Omega$ and
$$
\eta\omega=e^{Ku}(Nd\omega+\omega-\phi u_k \sigma^{kl}d_l)
$$
Since $\gamma= Dd$ on $\P\Omega$, $u_\gamma=\phi(x)\omega$ is equivalent to $\la Du, Dd\ra=\phi(x)\omega$ on $\P\Omega$. Differentiating $\la Du, Dd\ra=\phi(x)\omega$ with respect to $\P_1$ we find that
$$
\la Du_1, Dd\ra +\la Du, Dd_1\ra=\phi_1\omega+\phi\omega_1
$$
and so
\be \label{u1n}
u_{1n}(x_0)=\phi_1\omega+\phi\omega_1-u_1d_{11}
\ene
where we use the fact $\P_n=\gamma=Dd$ at $x_0$ and equation \eqref{der:d}.
Again by  $\eta=e^{Ku}(1-\phi^2)$ on $\P\Omega$, we get
\be \label{eta-1}
\eta_1=e^{Ku}((1-\phi^2)Ku_1-2\phi\phi_1)
\ene
On the other hand, the $\P_n$ derivative of $\omega$ at $x_0$ is
\be\label{omega_1n}
\omega_n=\F{u_1 u_{1n}}{\omega}+\F{u_nu_{nn}}{\omega}=\F{u_1}{\omega}u_{1n}+\phi u_{nn}
\ene
by assumption \eqref{u:assumption} and the fact $u_n=\phi\omega$ at $x_0$.
Combining identity \eqref{u1n} with equations \eqref{eta-1}, \eqref{omega1} and \eqref{omega_1n}, we obtain that
\be\label{fvomega}
\omega_{n}-\phi u_{nn}=\phi_1 u_1\F{1+\phi^2}{1-\phi^2}-\F{u_1^2}{\omega}d_{11}-u_1^2\phi K
\ene
We have $u_1^2=\omega^2(1-\phi^2)-1$ because of $u_n=\phi(x_0)\omega$ and assumption \eqref{u:assumption}. Combining these facts with equations \eqref{fvomega} and \eqref{comomega}, we find that
\begin{align*}
0 &\geq (K\phi\omega^2(1-\phi^2)+N\omega+\phi u_1\F{1+\phi^2}{1-\phi^2}-\F{u_1^2}{\omega}d_{11}- u_1^2\phi K -\phi_n\phi\omega)\\
&=\omega(\F{K\phi}{\omega}+N+\F{\phi u_1}{\omega}\F{1+\phi^2}{1-\phi^2}-\F{u_1^2}{\omega^2}d_{11}-\phi_n\phi)\\
&\geq \omega(N-C-\F{K}{\omega})
\end{align*}
where $C$ is a nonnegative constant depending on $d(x)$ and $\phi(x)$. Therefore $\omega(x_0)\leq K$ if $N\geq C+1$. The proof is complete.
\ep
\indent Next we establish some formulas about $\omega, v^k$. The Ricci curvature of $M$ will appear in our computation due to the commuting formula for covariant derivatives (Lemma \ref{lm:inchange}). \\
\indent
Two useful identities are
\begin{gather}
u^i u^j R_{ikjp}u^p=R(Du, \P_k, Du, Du)=0\quad \text{for any $k$}\label{vanishform}\\
\label{eq:misd}
v^k_i=\F{\sigma^{kl}u_{li}}{\omega}-\F{u^ku^l u_{li}}{\omega^3}=\F{g^{kl}u_{li}}{\omega}=h^k_i
\end{gather}
Recall that $\omega=\sqrt{1+|Du|^2}$ and $ v^k=\F{u^k}{\omega}$.
\bl \label{basic:lm}  For $v^k$ and $\omega$, we have
\begin{gather}
g^{ij}\omega_{ij}=(|A|^2+Ric(v_M, v_M))\omega+\la v_M, \nb(H\omega)\ra +\F{2}{\omega}g^{ik}\omega_i\omega_k\label{eq:omega}\\
g^{ij}v_{ij}^k=\F{g^{ij}g^{kl}u_{ijl}}{\omega}-\F{2}{\omega}g^{ij}v^k_i\omega_j-|A|^2v^k+Ric (g^{kl}\P_l, v_M)\label{eq:muk}
\end{gather}
Here $Ric$ denotes the Ricci curvature tensor of $M$ and $|A|^2$ is the second fundamental form of graph(u). Here $v_M$ is given by
\be\label{def:vn}
v_M =\F{Du}{\omega}
\ene
\el
\br
Since $\bar{\Omega}$ is compact, $Ric(v_M, v_M)$ is bounded below by a constant independent of $u(x)$.
\er
\bp Applying Lemma \ref{lm:inchange} to $g^{ij}u_{kij}$ we observe that
$$
\F{g^{ij}u_{kij}}{\omega}=\F{g^{ij}u_{ikj}}{\omega}=\F{g^{ij}u_{ijk}}{\omega}+\F{g^{ij}R_{jkip}u^p}{\omega}
$$
With equation \eqref{vanishform} we have
$$
\F{g^{ij}}{\omega}R_{jkip}u^p=\F{\sigma^{ij}}{\omega}R_{jkip}u^p-\F{u^i u^j}{\omega^2}R_{jkip}u^p
=Ric(\P_k, v_M);
$$
This gives
\be\label{eq}
\F{g^{ij}u_{kij}}{\omega}=\F{g^{ij}u_{ijk}}{\omega}+Ric(\P_k, v_M)
\ene
Now we compute $g^{ij}\omega_{ij}$ as follows.
\begin{align*}
g^{ij}\omega_{ij}&=g^{ij}(\F{u^k u_{ki}}{\omega})_j\\
&=\F{g^{ij}}{\omega}(\sigma^{kl}-\F{u^lu^k}{\omega^2})u_{lj}u_{ki}+\F{u^k}{\omega}g^{ij}u_{kij}\\
&=g^{ij}g^{kl}\F{u_{lj}u_{ki}}{\omega^2}\omega+\F{g^{ij}u^ku_{ijk}}{\omega}+\F{1}{\omega}Ric(Du, Du)\quad \text{by equation \eqref{eq}} \\
&=|A|^2\omega+\F{g^{ij}u^ku_{ijk}}{\omega}+Ric(v_M, v_M)\omega
\end{align*}
Then equation \eqref{eq:omega} follows from the following derivation
\begin{align*}
\la v_M,\nabla(H\omega)\ra &= \F{u^k}{\omega}(g^{ij}u_{ij})_k=g^{ij}\F{u^k}{\omega}u_{ijk}+\F{u^k}{\omega}\F{\P g_{ij}}{\P x_k}u_{ij}\\
&=g^{ij}\F{u^k}{\omega}u_{ijk}-\F{2}{\omega}g^{ik}\omega_i\omega_k
\end{align*}
The last   ne above is according to \eqref{parder:g}.
As for $g^{ij}v_{ij}^k$, equation  \eqref{eq:misd} yields that
\begin{align*}
g^{ij}v_{ij}^k &=g^{ij}(\F{g^{kl}u_{li}}{\omega})_j=g^{ij}(\F{g^{kl}}{\omega}u_{lij}-\F{g^{kl}}{\omega^3}u_{li}u^pu_{pj}+\F{1}{\omega}\F{\P g^{kl}}{\P u_p}u_{li}u_{pj}) \\
&=\F{g^{kl}}{\omega}g^{ij}u_{ijl}+Ric(g^{kl}\P_l, v_M)-\F{g^{ij}}{\omega}v^k_iw_j-\F{1}{\omega^3}g^{ij}(g^{kp}u^l+g^{lp}u^k)u_{li}u_{pj}\\
&=\F{g^{kl}}{\omega}g^{ij}u_{ijl}-2\F{g^{ij}}{\omega}v^k_i\omega_j-|A|^2v^k+Ric(g^{kl}\P_l, v_M)
\end{align*}
In the second line above we have applied equations \eqref{parder:g} and \eqref{eq}. The proof is complete.
\ep
\section{Estimates along the flows}
This section is devoted to establish two estimates about $u_t$ and $\omega$ along the flow in problem \eqref{eq:first}. First we derive the evolution equations of $u_t$ and $\omega$. Then we will explain how we follow the idea by \cite{GB96} mentioned in the introduction. Consequently the following two results are established.
\bt\label{thm:bound_ut} Let $u(x,t)$ be the smooth solution in problem \eqref{eq:first} on $\bar{\Omega}\TIMES [0,  T]$. Provided $\psi_{u}(x, u, Du)\geq -\chi_0$ on $\bar{\Omega}\TIMES [0,  T]$ for some $T>0$, then
\be
\max_{\bar{\Omega}\TIMES[0, T]}e^{-\chi_0 t}|u_t(x,t)|=\max_{x\in \bar{\Omega}}|u_t(x, 0)|
\ene
Here $\chi_0$ is a fixed constant.
\et
\br\label{Rm:ed} The case of $\chi_0=0$ in Euclidean spaces was obtained in (\cite{AW94},Lemma 2.2).
\er
\bt\label{thm:boundge}  Let $u(x,t)$ be the smooth solution of problem \eqref{eq:first} on $[0,  T]$ and $U_T$ denote $\max_{\bar{\Omega}\TIMES[0,T]}|u(x,t)|$. Assume $\psi(x, u, Du)$ is admissible with some constant $C$ (in page \pageref{condition:c1}). There exists a sufficiently large $K$ such that
\be
\omega\leq C_0e^{2K U_T}
\ene
on $\bar{\Omega}\PLH [0, T]$
where $C_0$ and $K$ are two constants only depending on the constant $C$, $d(x)$ and $\phi(x)$ and the Ricci curvature of $M$.
\et
\br $d(x)$ is defined in page \pageref{def:dx}.
\er
We define a parabolic operator
\be
L=g^{ij}\nb_i \nb_j-\psi_{u_k}\nb_k-\P_t
\ene
where $\nb$ is the covariant derivative of $M$. Comparing to Euclidean spaces, there are additional terms involving Ricci curvature in evolution equations of $\omega$ and $v^k$. See section 2 in \cite{GB96}.
\bl \label{eq:lm} Let $u(x,t)$ be the smooth solution to problem \eqref{eq:first} on $[0,  T]$. Then $u_t, \omega$ and $v^k$ satisfy that
\begin{gather}
Lu_t=\F{2}{\omega}g^{il}\omega_iu_{tl}+\psi_uu_t;\label{ev:ut}\\
L\omega=(|A|^2+(Ric(v_M, v_M)+\psi_u\la v_M,v_M\ra)\omega+\F{2}{\omega}g^{il}\omega_i\omega_l+\F{u^k}{\omega}\psi_k;\label{ev:omega}\\
Lv^k=(-|A|^2+\psi_u\F{1}{\omega^2})v^k+Ric(g^{kl}\P_l, v_M)+\F{g^{kl}}{\omega}\psi_l;\label{ev:vk}
\end{gather}
where $v_M$ is defined in Lemma \ref{basic:lm}.
\el
\bp A direct computation yields that
\begin{align*}
\P_t(u_t)&=\P_t(g^{ij}u_{ij})-\P_t(\psi(x,u,Du)); \\
&=g^{ij}(u_t)_{ij}-\F{2}{\omega}g^{il}\omega_i (u_t)_l-\psi_u u_t-\psi_{u_i}(u_t)_i \quad\text{by  \eqref{parder:g}}
\end{align*}
Equation \eqref{ev:ut} follows from reorganizing terms in the equation above. As for $\P_t\omega$ we find that
\begin{align*}
\P_t\omega&=\la v_M, \nb(H\omega)\ra-\F{u^k}{\omega}\psi_k-
\psi_u\la v_M,v_M\ra\omega-\psi_{u_i}\F{u_{ik}u^k}{\omega}\\
&=\la v_M,\nb(H\omega)\ra-\F{u^k}{\omega}\psi_k-
\psi_u\la v_M,v_M\ra\omega-\psi_{u_i}\omega_i
\end{align*}
With equation \eqref{eq:omega} this implies equation \eqref{ev:omega}. By $v^{k}=\F{u^k}{\omega}$ we derive
\begin{align}
\P_t v^k &=\F{\sigma^{kl}}{\omega}(\P_t u)_l-\F{u^k u^l}{\omega^3}(\P_t u)_l=\F{g^{kl}}{\omega}(\P_tu)_l\notag\\
&=\F{g^{kl}}{\omega}(g^{ij}u_{ij})_l- \F{g^{kl}}{\omega}\psi_l-\psi_u\F{g^{kl}u_l}{\omega}-\F{g^{kl}u_{il}}{\omega}\psi_{u_i}\label{vk:ms1}
\end{align}
Note that $g^{kl}u_l=(\sigma^{kl}-\F{u^ku^l}{\omega^2})u_l=\F{u^k}{\omega^2}$ and
\begin{align*}
\F{g^{kl}}{\omega}(g^{ij}u_{ij})_l &=\F{g^{kl}}{\omega}g^{ij}u_{ijl}-\F{2}{\omega}g^{ir}\omega_i u_{rl}\F{g^{kl}}{\omega}\\
&=\F{g^{kl}}{\omega}g^{ij}u_{ijl}-\F{2}{\omega}g^{il}\omega_i v^{k}_l
\end{align*}
Combining this with equation \eqref{vk:ms1} we obtain that
$$
\P_t v^k =\F{g^{kl}}{\omega}g^{ij}u_{ijl}-\F{2}{\omega}g^{il}\omega_i v^{k}_l-\F{\psi_u}{\omega^2} v^k-\psi_{u_i}v^k_i
-\F{g^{kl}}{\omega}\psi_l;
$$
Putting equation \eqref{eq:muk} and the above equation together we establish equation \eqref{ev:vk}. The proof is complete.
\ep
\indent Next we apply equation \eqref{ev:ut} to show Theorem \ref{thm:bound_ut}. The proof in Euclidean space shall still work in general Riemannian manifolds since there is no curvature term in equation \eqref{ev:ut}.
\bp Set $\vp=e^{-\chi_0 t}u_t$. By equation \eqref{ev:ut} $\vp$ satisfies that
\be\label{ev:vp}
L\vp=\F{2}{\omega}g^{il}\omega_i\vp_l+(\psi_u+\chi_0)\vp;
\ene
Fix $T>0$. In order to obtain Theorem \ref{thm:bound_ut} it suffices to prove if
$$
\vp(x_0, t_0)=\max_{\bar{\Omega}\TIMES[0,T]}\vp\geq 0
$$
then $t_0=0$. Suppose $t_0>0$. By our assumptions $\psi_u+\chi_0\geq0$, the maximal principle of parabolic equations implies that $x_0\in \P\Omega$. Now we take the coordinate $\{x_i\}$ in Remark \ref{cordinate}. Since $\vp$ obtains the maximum on $\bar{\Omega}$ at $(x_0,t_0)$, then
\be
\vp_k=u_{tk}=0\quad \text{for $1\leq k\leq n-1$}
\ene
Notice that $u_\gamma=\phi(x)\omega$ implies that $\la Du, Dd\ra =\phi \omega$ on $\P\Omega$. Taking $t$-derivative yields that
$$
u_{tn}=\phi(x_0)\omega_t
$$
On the other hand,
$$
\omega_t=\F{(\sum_{k\leq n-1}u_k u_{tk})+u_nu_{tn}}{\omega}=\phi(x_0)u_{tn}
$$
where we apply $u_{tn}=u_{nt}$ since $[\P_t,\P_n]=0$.
Hence $\omega_t=\phi^2(x_0)\omega_t$. Since $|\phi(x_0)|\leq \phi_0 <1$ $\omega_t=u_{tn}=0$. Thus $\vp_n(x_0, t_0)=0$. By Remark \ref{cordinate}, $\P_n$ at $x_0$ is the interior normal to $\P\Omega$. This is a contradiction to the Hopf Lemma unless $\vp$ is a constant on $\bar{\Omega}\TIMES [0, T]$. But the latter case implies that $t_0=0$. The proof is complete.
\ep
Our proof of Theorem \ref{thm:boundge} is based on the following estimate. This is the essence of the method in (\cite{GB96}).
\bl\label{lm:bounda} Suppose $\psi(x,u,Du)$ is admissible for some constant $C\geq 0$ ( assumptions (c1) in page \pageref{condition:c1}). Let $u(x,t)$ be the solution of problem \eqref{eq:first}. Define
$$
\eta =e^{Ku(x,t)}(Nd(x)+1-\phi(x)v^k d_k)
$$
where $N$ is from Lemma \ref{lemma:N}. There exists a sufficiently large $K$ with the property
if
$$(\omega\eta)(x_0,t_0)=\max_{\bar{\Omega}\TIMES[0,T]}(\omega\eta)
$$
for $x_0\in \Omega$ and fixed $T>0$,
then $\omega(x_0, t_0)\leq C_0$ for some constant $C_0$. Here $K$ and $C_0$ only depend on $C$, $d(x)$ and $\phi(x)$ and the Ricci curvature of $M$.
\el
This result generalizes Lemma 2.3 of \cite{GB96}.
\br\label{rm:onlyatpoint} A critical trick is that Ricci curvature on $\bar{\Omega}$ is absorbed in the process to  determine $K$. See equations \eqref{final:c} and \eqref{Final:c}.
\er
\bp Fix any $T>0$. Now assume $x_0\in\Omega$ and $\eta\omega$ achieves its maximum  on the set $\{ (x,t):\bar{\Omega}\TIMES [0, T]\}$ at $(x_0, t_0)$ for $t_0>0$. Again all expressions are evaluated at $(x_0, t_0)$. In the following we denote different positive constants by $C_1$, having no dependence on $K$. \\
\indent First we observe that  $\tw_i\eta+\tw\eta_i=0$ for $i=1,\cdots,n$ and
\be
0\geq L(\tw\eta)=(L\tw-\F{2}{\tw}g^{il}\tw_i\tw_l)\eta+\tw L\eta
\ene
Combining inequality \eqref{eq:csd} with equation \eqref{ev:omega} and dividing it by $\eta\omega$ we find that
\be\label{eq:csd}
\F{L\eta}{\eta}+|A|^2+Ric(v_M, v_M)+\psi_u\la v_M, v_M\ra+\F{u^k}{\omega^2}\psi_k\leq 0
\ene
By assumptions \eqref{condition:c1} the Cauchy inequality gives
$$|u^k\psi_k|\leq \omega |\psi_k|\leq C\omega^2, \psi_u\geq -C $$ Since $\bar{\Omega}$ is compact and $|v_M|\leq 1$,  we get
\be\label{final:c}
Ric(v_M, v_M)+\F{u^k}{\omega^2}\psi_k+\psi_u\la v_M,v_M\ra\geq -C_1
\ene
for a constant $C_1$ determined by $C$. Consequently inequality \eqref{eq:csd} may be simplified further into
\be\label{Final:c}
0\geq  \F{L\eta}{\eta}+|A|^2-C_1
\ene
In order to analyze $L\eta$ we write $\eta$ as $e^{Ku}h$ where $h=Nd+1-\phi v^kd_k$. If denote $\max_{x\in \bar{\Omega}}\{d(x)\}$ by $d_0$, $h$ satisfies that
\be\label{bound:h}
1-\phi_0\leq h\leq (Nd_0+1+\phi_0)
\ene
Thus $\F{L\eta}{\eta}$ may be expanded as follows.
\begin{gather}\label{eta:term}
\F{L\eta}{\eta}=\F{K^2}{h}g^{ij}u_iu_j+KLu+2K \F{g^{ij}u_i h_j}{h}+\F{Lh}{h}
\end{gather}
We can assume $ \omega(x_0, t_0)\geq \sqrt{2}$. Otherwise we are done. Thus
$$
g^{ij}u_iu_j=(\sigma^{ij}-\F{u^i u^j}{\omega^2})u_iu_j=\F{|Du|^2}{\omega^2}\geq \F{1}{2}
$$
With inequality \eqref{bound:h} the first term in equation \eqref{eta:term} is bounded below as follows.
\be\label{computing:A}
\F{K^2}{h}g^{ij}u_iu_j\geq \F{K^2}{(Nd_0+1+\phi_0)}\F{|Du|^2}{\omega^2}\geq C_2 K^2
\ene
where $C_2:=(2(Nd_0+1+\phi_0))^{-1}$. By assumptions \eqref{condition:c1}, the second term in equation \eqref{eta:term} becomes
\be\label{computing:B}
K Lu=K(\psi(x,u, Du)-\psi_{u_i}(x, u, Du)u_i)\geq -CK
\ene
The following identities about $h$ are easily verified.
\begin{gather}
\P_t  h=-\phi(x)v^k_t d_k\label{eq:h:t}\\
h_i =Nd_i-(\phi d_k)_i v^k -\phi d_k v^k_i\label{eq:h:i}\\
h_{ij} = N d_{ij}-(\phi d_k)_{ij}v^k-(\phi d_k)_i v^k_j-(\phi d_k)_j v^k_i -\phi d_k v_{ij}^k\label{eq:h:ij}
\end{gather}
A  direct computation yields that $g^{ij}u_j=\F{u^i}{\omega^2}$.   Applying the Cauchy inequalities \eqref{Cauchy:Inequality} and \eqref{estimate:G}, with equation \eqref{eq:h:i}, the third term in equation \eqref{eta:term} becomes
\begin{align}
2K\F{g^{ij}u_i h_j}{h}&=2\F{K}{h\omega^2}(u^iNd_i-u^i(\phi d_k)_i v^k -u^i\phi d_k v^k_i)\notag\\
&\geq -\F{C_1K}{\omega}-\F{C_1K|A|}{\omega}\notag\\
&\geq -C_1 K-\F{C_1K^2}{4\Sc\omega^2}-\Sc|A|^2\label{computing:C}
\end{align}
where $C_1$ denotes some constant independent of $K$. Here $\Sc$ is a small constant determined later. Here we also used the fact $|v^k_i|=|A|$.\\
\indent Next we consider the last term in equation \eqref{eta:term}. Substituting $Lv^k$ with equation \eqref{ev:vk} we find that
\begin{align}
\F{Lh}{h}&=\F{1}{h}\{g^{ij}(Nd_{ij}-(\phi d_k)_{ij}v^k-(\phi d_k)_i v^k_j-(\phi d)_j v^k_i)-\psi_{u_i}(Nd_i-(\phi d_k)_iv^k) -\phi d_k Lv^k\}\notag\\
&\geq -C_1-C_1|A|+|A|^2\F{\phi d_k}{h}v^k-\F{\phi d_k}{h}Ric(g^{kl}\P_l, v_M)\notag\\
&\geq -C_1(2+\F{1}{4\Sc})-\Sc|A|^2+\F{\phi d_kv^k}{h}|A|^2\label{computing:D}
\end{align}
The second inequality follows from $|v^i_k|=|A|$ and the Cauchy inequalities in \eqref{estimate:G}. As for the third line, notice that
$$|\F{\phi d_k }{h}g^{kl}Ric(\P_l, v_M)|\leq C_1$$
by assumptions \eqref{condition:c1} and the compactness of $\bar{\Omega}$. Applying inequality \eqref{estimate:G} and combining these facts together we will obtain the third line. \\
\indent Finally we put the estimates in inequalities \eqref{computing:A}, \eqref{computing:B}, \eqref{computing:C} and \eqref{computing:D} into equation \eqref{eta:term} and choose $C_1>C\geq 0$. Then we conclude that
\begin{align}
\F{L\eta}{\eta}&\geq K^2(C_2-\F{C_1}{4\Sc\omega^2})-2C_1 K-C_1(2+\F{1}{4\Sc})+(\F{\phi d_kv^k}{h}-2\Sc)|A|^2\label{eq:uv}
\end{align}
Now we choose $2\Sc=\F{1}{Nd_0+1+\phi_0}$. Thus
$$
\F{Nd+1}{Nd+1-\phi v^k d_k}-2\Sc\geq 0
$$
With inequality \eqref{eq:uv}, inequality \eqref{Final:c} leads that
\begin{align*}
0 &\geq  K^2(C_2-\F{C_1}{4\Sc\omega^2})-2C_1 K-C_1(2+\F{1}{4\Sc})+|A|^2(\F{Nd+1}{Nd+1-\phi v^k d_k}-2\Sc)\\
&\geq K^2(C_2-\F{C_1}{4\Sc\omega^2})-2C_1 K-C_1(3+\F{1}{4\Sc})
\end{align*}
Taking $K$ sufficiently large, we obtain $\omega(x_0,t_0)\leq C_0$ where $C_0$ relies on $K, N$, $d(x)$ and $\phi(x)$. Moreover $K$ is determined by $N$, $C$ and the Ricci curvature of $M$.  As pointed by Lemma \ref{lemma:N}, $N$ depends on $d(x)$ and $\phi(x)$. Thus the proof is complete.
\ep

The preceding result allows us to conclude Theorem \ref{thm:boundge}.
\bp(The proof of Theorem \ref{thm:boundge}) Suppose N is from Lemma \ref{lemma:N} and $K$ is given by Lemma \ref{lm:bounda}. We consider the maximum of $\omega\eta$. Fix $T>0$, assume
$$
\omega\eta (x_0, t_0)=\max_{x\in\bar{\Omega}, t\in [0, T]}(\omega\eta)
$$
If $x_0\in \P\Omega$, Lemma \ref{lemma:N} says that $\omega(x_0,t_0)\leq K$. Otherwise, Lemma \ref{lm:bounda} indicates that $\omega(x_0, t_0)\leq C_0$ for $x_0\in \Omega$. Now assume $C_0\geq K$. We have $\omega(x_0, t_0)\leq C_0$ in both cases.
Thus for any point $(x,t)$ in $\bar{\Omega}\times [0, T]$
$$
\omega \leq e^{K(u(x_0,t_0)-u(x,t))}\omega(x_0,t_0)\leq C_0e^{2KU_T}
$$
where $U_T$ denotes the maximum of $|u(x,t)|$ on $\bar{\Omega}\times [0, T]$. From Lemma \ref{lm:bounda} $K$ and $C_0$ only rely on the admissible constant $C$, $\phi(x)$, $d(x)$ and the Ricci curvature of $M$. Hence we conclude Theorem \ref{thm:boundge}.
\ep

\section{Convergence of the flows}
In this section we discuss long time existences and asymptotic behaviors of solutions to problem \eqref{eq:first} via the estimates in the previous section. We will establish Theorem \ref{thm:GLT}, Theorem \ref{thm:MT1} and Theorem \ref{thm:MT2}.\\
\indent  A sufficient condition for long time existences of solutions to problem \eqref{eq:first} is that problem \eqref{eq:first} preserves strictly parabolic on any finite time interval $[0, T]$. Naturally our goal is to establish
$
\omega\leq C
$
on $\bar{\Omega}\TIMES [0, T)$ where $C$ is a finite constant possibly depending on $T$.
This idea is realized in the following result, as a generalization of Theorem 2.4 in \cite{GB96}.
\begin{customthm}{\ref{thm:GLT}}  Suppose $\psi(x,u,Du)$ is a smooth function of $x,u$ and $Du$ satisfying 
	\begin{enumerate}
		\item [(a)]$\psi_u(x,u,Du)\geq c_0$ where $c_0$ does not rely on $u$;
		\item [(b)]for any given $K>0$ there is a positive constant $C=C(K)$ such that $\psi(x,u,Du)$ is admissible with respect to $C$ if $|u|\leq K$.
	\end{enumerate}  Then problem \eqref{eq:first} with $\psi(x,u,Du)$ has a solution $u\in C^{\infty}(\bar{\Omega}\times [0,\infty))$.																																										
\end{customthm}
\br Two examples about admissible functions are given in Lemma \ref{lm:condition}.
\er
\bp The short time existence of the solution in problem \eqref{eq:first} follows from the standard theory of parabolic equations. Let $T<\infty$ be the maximal time such that the solution $u(x,t)$ of problem \eqref{eq:first} exists smoothly on $[0, T)$. By Theorem \ref{thm:bound_ut} condition (a) implies that 
$$|u_t|\leq C_1 e^{c_0 T}$$
for all $(x,t)\in \bar{\Omega}\times [0, T)$. Here $C_1$ is a constant depending on the initial condition $u_0(x)$. Let $U$ denote the supremum of $|u(x,t)|$ on $\bar{\Omega}\times [0, T)$.
Thus
$$
U\leq \max_{x\in \bar{\Omega}, t\in [0, T)}|u(x,t)-u(x,0)|+\max_{x\in\bar{\Omega}}|u(x,0)|\leq C_1 e^{c_0 T}T+U_0
$$
where $U_0=\max_{x\in\bar{\Omega}}|u(x,0)|$. Since $\psi(x,u,Du)$ is admissible on $[0, T)$ for any $u(x,t)$, Theorem \ref{thm:boundge} gives the gradient estimate 
$$
\omega\leq C_2:=C_0 e^{2K(C_1 e^{c_0T}T+U_0)}
$$
Thus problem \eqref{eq:first} is strictly parabolic until time $T$. Thus $u(x,t)$ can be extended over time $T$. This gives a contradiction. The proof is complete.
\ep
Letting $\psi\equiv 0$ we obtain long time existence of mean curvature flows of graphs with prescribed contact angle in product manifolds.
\begin{cor}\label{cor:LT}
	Problem \eqref{eq:first} with $\psi\equiv 0$ has a solution $u(x,t)\in C^{\infty}(\bar{\Omega}\TIMES [0,\infty))$.
\end{cor}
We give a sufficient condition for uniformly convergence of the solution to problem \eqref{eq:first} as $t\rightarrow \infty$.
\bl\label{lm:cn}  Let $u\in C^{\infty}(\bar{\Omega}\times [0,\infty))$ be the solution to problem \eqref{eq:first}. If
\begin{itemize}
	\item [(i)]$u(x,t)$ and its higher derivatives are uniformly bounded;
	\item [(ii)]for any $T>0$, $\int_{0}^T\int_{\bar{\Omega}}|u_{t}|^2 dxdt\leq C$ for some uniformly constant $C$;
\end{itemize}
then $u_t(x,t)$ converges uniformly to $0$ as $t\rightarrow \infty$ and
\begin{itemize}
	\item [(i)] for $\psi(x,u,Du)\equiv 0$, $u(x,t)$ converges uniformly to a constant,
	\item [(ii)]for $\psi(x,u, Du)\equiv h(x,u)\omega$ with $h_u(x,u)\geq h_0>0$, $u(x,t)$ converges uniformly to a smooth function $u_{\infty}(x)$ which is the solution to problem \eqref{eq:second}.
\end{itemize}
\el
Its proof is presented in appendix D.\\
\indent Now we extend the result of \cite{Hui89A} in Euclidean space into general Riemannian manifolds. The uniformly convergence is still valid without any restriction on Ricci curvature.
\begin{customthm}{\ref{thm:MT1}} Problem \eqref{eq:first} with $\psi\equiv \phi\equiv 0$ has a smooth solution $u(x,t)$ on $\bar{\Omega}\times[0, \infty)$. Moreover $u(x,t)$ converges uniformly to a constant as $t$ goes to infinity.
\end{customthm}
In this case problem \eqref{eq:first} becomes
\be\label{st:max}
\left\{\begin{aligned}
	&u_t =g^{ij}(Du)u_{ij},\quad \text{on $\Omega\times [0,\infty)$}\\
	&u_{\gamma} =0\quad (x,t)\in \P\Omega\times [0,\infty), \quad  u(x,0)=u_0(x)\quad x\in\Omega\\
\end{aligned}
\right.
\ene
where $g^{ij}=(\sigma^{ij}-u^iu^j\omega^{-2})$.
It describes a nonparametric mean curvature flow with a vertical contact angle.
First the $C^0$ bound of the solution to problem \eqref{eq:first} is an immediate result about the strong maximum principle (see \cite{Lie96}).
\bl\label{lm:control:u}  Let $u(x,t)$ be the smooth solution in \eqref{eq:first} with $\psi\equiv \phi\equiv 0$. Then
\be\label{eq:control:u}
\min_{x\in \bar{\Omega}}u_0(x)\leq u(x,t)\leq \max_{x\in \bar{\Omega}}u_0(x);
\ene
for any $t$.
\el
\bp Fix any $T>0$. It suffices to show that if
$$ u(x_0, t_0)=\max_{(x,t)\in \Omega \times[0,T]}u(x,t)\geq 0 $$
then $t_0=0$. Suppose $t_0>0$. The maximal principle, with equation \eqref{st:max}, implies $x_0\in \P\Omega$. However, $u_\gamma(x_0)=0$. The Hopf Lemma or the strong maximal principle (see \cite{Lie96}) implies that $u(x,t)$ is a constant. This means $t_0=0$ and gives a contradiction. Therefore $t_0=0$. The proof is complete.
\ep
Now we investigate $u_t$ along the flows. For later reference, we work with much general $\psi(x,u,Du)$.
\bl \label{uni:k}  Suppose $\psi(x,u,Du)$ is $h(x,u)\omega$. Let $u(x,t)$ be the smooth solution to problem \eqref{eq:first} on $[0, T]$ with $u_\gamma =\phi(x)\omega$ on $\P\Omega$. Then
\be
\int_0^T\int_{\bar{\Omega}}\F{(u_t)^2}{\omega}dxdt=-
\{\int_{\bar{\Omega}}\omega dx+\int_{\P\Omega}u\phi(x)dx+\int_{\bar{\Omega}}g(x,u)dx\}|_0^T
\ene
where $g(x,u)=\int h(x,u)du$.
\el
\bp In this setting, problem \eqref{eq:first} is rewritten as
$$
u_t=div(\F{Du}{\omega})\omega-h(x,u)\omega
$$
Thus the divergence theorem yields that
\be\label{div:f1}
div(u_t\F{Du}{\omega})=\la \F{Du}{\omega}, D(u_t)\ra +\F{u_t}{\omega}u_t+h(x,u)u_t;
\ene
Next we compute the $t$-derivative of $\int_{\bar{\Omega}}\omega dx$.
\begin{align*}
\P_t \int_{\bar{\Omega}}\omega dx&=\int_{\bar{\Omega}}\la \F{Du}{\omega}, D(u_t)\ra dx\\
&= \int_{\bar{\Omega}}div(u_t\F{Du}{\omega})dx -\int_{\bar{\Omega}}\F{u_t}{\omega}(u_t+h(x,u)\omega)dx\\
&=-\int_{\P\Omega}u_t\phi(x)dx-\int_{\bar{\Omega}}\F{(u_t)^2}{\omega}dx-\P_t\int_{\bar{\Omega}}g(x,u)dx
\end{align*}
where $g(x,u)=\int h(x,u) du$. Reorganizing the above equation, we complete the proof.
\ep
Now we are ready to conclude Theorem \ref{thm:MT1}.
\bp
By Lemma \ref{lm:control:u}, we get a uniform bound of $u$. Notice that $\psi\equiv 0$ is admissible with the constant $0$. Consequently $u_t$ and $\omega$ are uniformly bounded by Theorem \ref{thm:bound_ut} and Theorem \ref{thm:boundge} respectively. Thus $ u\in C^{\infty}(\bar{\Omega}\TIMES [0,\infty))$. As a result by problem \eqref{eq:first} all derivatives of $u(x,t)$ are uniformly bounded. From Lemma \ref{uni:k} and uniform bounds of $u$ and $\omega$, we observe that
$$
\int_0^T\int_{\bar{\Omega}}|u_t|^2 dx dt\leq C
$$
Thus all conditions of Lemma \ref{lm:cn} are fulfilled. Consequently $u(x,t)$ converges uniformly to a constant by Lemma \ref{lm:cn}.
\ep
Now we consider a little complicated form of $\psi(x,u,Du)$, i.e $\psi(x,u, Du)=h(x,u)\omega$. 
\begin{customthm}{\ref{thm:MT2}}  Let $\psi(x,u, Du)$ be $h(x,u)\omega$ with $h_u\geq h_0>0$. Problem \eqref{eq:first} has a smooth solution $u(x,t)$ on $\bar{\Omega}\times [0,\infty)$. Moreover $u(x,t)$ converges uniformly to $u_\infty(x)$ as the solution to problem \eqref{eq:second} as $t$ goes to infinity.
\end{customthm}
Recall that $g^{ij}(Du)=(\sigma^{ij}-u^iu^j\omega^{-2})$.
In this case, problem \eqref{eq:first} takes the form
\be
\left\{\begin{aligned}
	&u_t =g^{ij}(Du)u_{ij}-h(x,u)\omega,\quad\text{on $\Omega\times [0,\infty)$}\\
	&u_{\gamma} =\phi(x)\omega,\quad  \quad \text{in $ \P\Omega\TIMES \ [0,\infty)$}\\
	&u(x,0)=u_0(x),\quad x\in \Omega
\end{aligned}
\right.
\ene
and problem \eqref{eq:second} is rewritten as
\be
g^{ij}(Du)u_{ij}=h(x,u)\omega \quad\text{on $\Omega$}
\ene
with $u_\gamma=\phi(x)\omega$ on $\P\Omega$.
First we establish the $C^0$ bound with a little more general assumption. The following result shows that positive gravity is a very strong condition.
\bl\label{lm:est:u}
Suppose $u(x,t)$ is the solution to problem \eqref{eq:first} with $\psi_u\geq \psi_0>0$. Then
\begin{itemize}
	\item [(i)] $|u(x,t)|\leq \max_{x\in \bar{\Omega}}|u(x,0)|+\F{1}{\psi_0}\max_{x\in\bar{\Omega}}|u_t(x,0)|$.
	\item [(ii)]Fix $m>0$. For any $\psi_0\in (0, m]$  $$\psi_0|u(x,t)|\leq m\max_{x\in \bar{\Omega}}|u(x,0)|+\max_{x\in\bar{\Omega}}|u_t(x,0)|$$
\end{itemize}
\el
\bp Fix any $t>0$. Since $\psi_u\geq \psi_0>0$, Theorem \ref{thm:bound_ut} implies that
\be\label{ut:decay}
|u_t|\leq e^{-\psi_0 t}\max_{x\in\bar{\Omega}}|u_t(x,0)|
\ene
Thus
$$
|u(x,t)|\leq \max_{x\in \bar{\Omega}}|u(x,0)|+ t e^{-\psi_0 t}\max_{x\in\bar{\Omega}}|u_t(x,0)|
$$
The conclusion (i) follows from $te^{-\psi_0 t}\leq \F{1}{\psi_0}$.  Multiplying both sides by $\psi_0$, we obtain the conclusion (ii).
\ep
Now it is time to show Theorem \ref{thm:MT2}.
\bp Let $u(x,t)$ be the solution to problem \eqref{eq:first}. Since $h_u\geq h_0>0$, $\psi_u=h_u(x,u)\omega>h_0>0$. Thus Lemma \ref{lm:est:u} implies that $u(x,t)$ is uniformly bounded. We obtain that $\psi(x,u, Du)$ is admissible for a fixed constant by Lemma \ref{lm:condition}. Consequently Theorem \ref{thm:boundge} shows that $\omega$ is uniformly bounded. Hence $u(x,t)\in C^\infty(\bar{\Omega}\TIMES[0,\infty))$. Since $\psi_u\geq 0$ we observe that $u_t$ are uniformly bounded by Theorem \ref{thm:bound_ut}. According to problem \eqref{eq:first} all high derivatives of $u(x,t)$ are uniformly bounded. From Lemma \ref{uni:k} these facts in turn imply
$$
\int_0^T\int_{\bar{\Omega}}|u_t|^2dxdt\leq C
$$
for any $T>0$.
By Lemma \ref{lm:cn} $u(x,t)$ converges uniformly to $u_{\infty}(x)$. Moreover $u_t$ converges to $0$. Thus $u_{\infty}(x)$ is a smooth solution to problem \eqref{eq:second}.
\ep
From viewpoint of partial differential equations we solve the Capillary problem with positive gravity via a flow method. Restating the convergence part of Theorem \ref{thm:MT2} gives that
\bt\label{thm:MT22} Let $\psi(x,u, Du)$ be $h(x,u)\omega$ with $h_u\geq h_0>0$. Problem \eqref{eq:second} has a unique smooth solution.
\et
\bp  From Lemma \ref{lm:unique} below we have the uniqueness. The existence part follows from Theorem \ref{thm:MT2}. We complete the proof.
\ep
Now we show the uniqueness in Theorem \ref{thm:MT22}.
\bl\label{lm:unique} Let $\psi(x,u, Du)$ be $h(x,u)\omega$ with $h_u\geq h_0>0$. There is at most one smooth solution to problem \eqref{eq:second}.
\el
\bp Notice that $\psi_u(x,u,Du)\geq h_0 >0$. Suppose $u_1(x)$ and $u_2(x)$ are two smooth solutions to problem \eqref{eq:second}. Consider $u(x)=u_1(x)-u_2(x)$. By Lemma \ref{difference} in appendix C $u(x)$ satisfies a quasilinear equation
$$
g^{ij}(Du_1) u_{ij}+\tilde{b}_iu_i=\psi_u(x,*, Du_1)u
$$
where $(g^{ij}(Du_1))$ is a positive definite matrix and $*$ denotes some smooth function. Since $\psi_u \geq h_0>0$, $u(x)$ can not be a constant function except $0$. By the maximal principle, the positive maximum of $u$ occurs at some $x_0\in \P\Omega$ or $u\leq 0$. For the first case, Lemma \ref{gamma:vanish} in appendix C implies that $u_\gamma(x_0)=0$. Thus the Hopf Lemma gives a contradiction since $u(x)$ can not be a positive constant function. Hence we obtain that $u\leq 0$. Reversing the role of $u_1$ and $u_2$ yields that $u$ has to be equal to zero. Therefore the solution to problem \eqref{eq:second} is unique.
\ep
Another application of Theorem \ref{thm:MT2} is stated as follows.
\begin{cor}\label{cor:est:u} Suppose $\psi(x,u,Du)$ is $\Sc u$. Let $u_\Sc(x)\in C^{\infty}(\bar{\Omega})$ be the solution to problem \eqref{eq:second}. If $\Sc\in (0, 1]$, then
	\be
	|\Sc u_{\Sc}(x)|\leq C_1
	\ene
	where $C_1$ is a constant independent of $\Sc$.
\end{cor}
\bp Suppose $u_0(x)$ is a smooth function on $\bar{\Omega}$ with $u_\gamma=\phi(x)\omega$ on $\P\Omega$. Let $u_{\Sc}(x,t)$ be the solution to problem \eqref{eq:first} with $\psi(x,u,Du)=\Sc u$ and the initial condition $u(0,x)=u_0(x)$. For $\Sc\in (0,1]$, from conclusion (2) in Lemma \ref{lm:est:u} we observe that
$$
\Sc |u_{\Sc}(x,t)|\leq C_1
$$
where $C_1$ only depends on $u_0(x)$. According to Theorem \ref{thm:MT2} $u_{\Sc}(x,t)$ converges uniformly to $u_{\Sc}(x)$ as $t\rightarrow \infty$. As a result $|\Sc u_{\Sc}(x)|\leq C_1$. The proof is complete.
\ep

\section{Translating Surfaces}
In the evolution of problem \eqref{eq:first} in Section 4, a common feature is that $u(x,t)$ is uniformly bounded. However not all solutions to problem \eqref{eq:first} have this feature even if $\psi(x,u,Du)$ is $0$. In this section we will  work with one such example.\\
\indent We consider problem \eqref{eq:first} with $\psi(x,u, Du)\equiv 0$ and $\Omega$ as a convex domain in a Riemannian surface $M^2$. We rewrite problem \eqref{eq:first} as
\be\label{eq:2:ts}
\left\{\begin{split}
	&u_t =g^{ij}(Du)u_{ij}\quad \text{on $\Omega\times [0,\infty)$}\\
	&u_{\gamma} =\phi(x)\omega,\quad  \quad (x,t)\in \P\Omega\TIMES \ [0,\infty)\quad
	u(x,0)=u_0(x),\quad x\in \Omega  \end{split}\right.
\ene
where $g^{ij}=(\sigma^{ij}-u^iu^j\omega^{-2})$. Throughout this section problem \eqref{eq:2:ts} is problem \eqref{eq:first} with $\psi\equiv 0$. \\
\indent The purpose of this section is to generalize the result of (\cite{AW94}) into nonegative Ricci curvature case as follows.
\begin{customthm}{\ref{thm:MT3}} Suppose $M^2$ is a Riemannian surface with nonnegative Ricci curvature. Let $\Omega$ be a bounded convex domain with
	$$k_0 \geq k\geq |\F{\phi_T}{\sqrt{1-\phi^2}}|+ \delta_0$$
	where $k$ is the inward curvature of $\P\Omega$ and $T$ is the tangent vector of $\P\Omega$, $k_0,\delta_0$ are positive constants.
	Then
	\begin{itemize}
		\item [(i)]problem \eqref{eq:first} with $\psi(x,u,Du)\equiv 0$ has a smooth solution $u(x,t)$ on $\bar{\Omega}\times[0,\infty)$ with  $\omega\leq C_1$ where $C_1$ is a constant depending on $ k_0, \delta_0, \phi_0$ and the initial condition.
		\item [(ii)]
		moreover $u(x,t)$ converges uniformly to $u_{\infty}(x)+Ct$ as $t\rightarrow \infty$ where $u_{\infty}(x)$ is the solution to problem \eqref{eq:second} with $\psi(x,u, Du)\equiv C$. Here $C$ is given by
		$$
		C=-\F{\int_{\P\Omega}\phi(x)dx}{\int_\Omega (1+|Du_{\infty}|^2)^{-\F{1}{2}}dx}
		$$
	\end{itemize}
\end{customthm}
In view of Corollary \ref{cor:LT} we already have the long time existence. Since conclusion (ii) above says that $u(x,t)$ is unbounded, Theorem \ref{thm:MT3} is totally different with the convergence results in Section 4.\\
\indent First we give some preliminary facts. Because $\P\Omega$ is a smooth curve, we need more special notation here. The following constructions always work in a neighborhood of a fixed point in $\P\Omega$. Without loss of generality we always assume them valid on $\P\Omega$.\\
\indent  Let $T$ be a unit smooth tangent vector of $\P\Omega$ and $\gamma$ still denote the interior normal to $\Omega$. The curvature of $\P\Omega$ is
$$
k=\la \nb_{T}T,\gamma\ra
$$
where $\nb$ is the covariant derivative of $M^2$. Let $\{\theta\}$ denote a local coordinate on $\P\Omega$ and $r(x)$ be the distance function $d(x,\P\Omega)$. As in Lemma \ref{eq:local coordinate}  we construct a frame $\{\P_r, \P_\theta\}$ in a neighborhood of $\P\Omega$ such that
$$\la \P_r, \P_\theta\ra=0,\quad \la \P_r,\P_r\ra=1$$
We denote the function $\la \P_\theta, \P_\theta\ra^{\F{1}{2}}$ by $\vp(x)$. Thus we obtain an orthonormal frame $\{\P_r, \vp^{-1}\P_\theta\}$ near $\P\Omega$ noted as $\{\gamma, T\}$ for short.\\
\indent Given a smooth function $f$ we write $f_X$ for $X(f)$ and $f_{XY}$ for the covariant derivative of $f$ defined by
$$
f_{XY}=X(Y(f))-(\nb_XY)f;
$$
where $X$ and $Y$ are tangent vector fields.
The following result is classical about the geometry of curves (see \cite{AW94}, \cite{Gra89}). For completeness we present its proof here.
\bl\label{lm:estd} Let $\Omega$ be a smooth domain in a Riemannian surface $M^2$ with covariant derivative $\nb$. Then on $\P\Omega$,
\begin{itemize}
	\item [(i)]$\nb_TT=k\gamma$, $\nb_T\gamma =-kT$, $\nb_\gamma T=\nb_\gamma\gamma=0$;
	\item [(ii)] for any $f\in C^{\infty}(\bar{\Omega})$,
	$
	\gamma(T(f))=T(\gamma(f))+kf_T
	$.
\end{itemize}
\el
\bp Note that $\la \nb_T T, T\ra \equiv 0$ on $\P\Omega$. The definition of $k$ yields that $\nb_TT=k\gamma$. And $\nb_T\gamma =-kT$ follows from $\la \nb_T\gamma,\gamma\ra=0$ and $\la \nb_T\gamma, T\ra=-k$. In addition we find
$$
\P_r\vp=\F{1}{\vp}\la \nb_{\P_r}\P_\theta, \P_\theta\ra=\vp\la \nb_{T}\gamma, T\ra=-k\vp
$$
Thus we conclude that
\begin{align*}
\nb_{\gamma}T &=\nb_{\P_r}\vp^{-1}\P_\theta=\vp^{-1}\nb_{\P_\theta}\P_r-\vp^{-2}\P_r\vp\P_\theta\\
&=\nb_T\gamma+kT=0
\end{align*}
By $\la\nb_{\gamma}T,\gamma\ra+\la T,\nb_{\gamma}\gamma\ra=0$, we have $\la T,\nb_{\gamma}\gamma\ra=0$. This gives $\nb_{\gamma}\gamma=0$ since $\la \gamma,\gamma \ra=1$. As for conclusion (ii), a direct computation yields that
\begin{align*}
\gamma(T(f))&=\P_r(\vp^{-1}\P_\theta f)=-\F{1}{\vp^2}\P_r\vp\P_{\theta} f+\vp^{-1}\P_\theta\P_r f\\
&=k\vp^{-1}\P_\theta f+T(\gamma (f))
\end{align*}
We establish conclusion (ii) and accomplish the proof.
\ep
Now we explain how the convex boundary of $\Omega$ controls $\omega$. In a neighborhood of $\P\Omega$, $\{T,\gamma\}$ is an orthonormal frame. Thus near $\P\Omega$ we can write
$$
|Du|^2=u^Tu_T+u^\gamma u_\gamma=u_T^2+u_\gamma^2
$$
Recall that $\omega=\sqrt{1+|Du|^2}$. In the following derivation all computations only happen on $\P\Omega$. From $u_\gamma =\phi(x)\omega$ we obtain
\be
u_\gamma^2=\F{\phi^2}{1-\phi^2}(u_T^2+1), \quad u_T^2=(1-\phi^2)\omega^2-1;
\ene
Consequently some second derivatives of $u$ (not covariant derivative) can be expressed as follows:
\be \label{ded}
\begin{split}
	TT(u)&=\F{-\phi\phi_T\omega^2+\omega\omega_T(1-\phi^2)}{u_T} \\
	T\gamma(u)&=\phi_T\omega+\phi\omega_T\\
	\gamma T(u)&=\phi_T\omega+\phi\omega_T+k u_T
\end{split}
\ene
Here we assume $u_T\neq 0$. The last identity follows from conclusion (ii) in Lemma \ref{lm:estd}. We can not compute $\gamma(\gamma(u))$ directly because we need more information away from the boundary. On $\P\Omega$, the expressions of $g^{\gamma T}$, $g^{TT}$ and $g^{\gamma\gamma}$ are given as follows.
\be
\begin{split}
	g^{TT}&=(1-\F{u^Tu^T}{\omega^2})=\F{1+u_\gamma^2}{\omega^2}\\
	g^{\gamma\gamma}&=(1-\F{u^\gamma u^\gamma }{\omega^2})=1-\phi^2(x)\\
	g^{\gamma T}&=g^{T\gamma}=-\F{u^\gamma u^T }{\omega^2}=-\F{u_\gamma u_T}{\omega^2}
\end{split}
\ene
The next result is essential in this section.
\bl\label{lm:convex:bound} Let $\Omega$ be a convex domain in $M^2$ with
\be
k_0\geq k\geq |\F{\phi_T}{\sqrt{1-\phi^2}}|+\delta_0
\ene
where $k$ is the inward curvature of $\P\Omega$, $k_0, \delta_0$ are positive constants.
Suppose $u\in C^\infty(\bar{\Omega})$ satisfies
\be
|g^{ij}u_{ij}|\leq C_0\quad u_\gamma=\phi(x)\omega,\quad |\phi|\leq \phi_0<1
\ene
for all $x\in \P\Omega$. If the maximum of $\omega$ occurs at $x_0\in \P\Omega$, then $\omega(x_0)\leq C_1$ where $C_1$ is a constant depending on $C_0$, $k_0,\delta_0$ and $\phi_0$.
\el
\br This phenomenon was firstly observed by \cite{AW94} in Euclidean spaces. Here we do not impose any Ricci curvature condition on $M$.
\er
\bp Suppose $\omega$ achieves its maximum at $x_0\in \P\Omega$. In the following all expressions are evaluated at $x_0$. Recall that on $\P\Omega$,
$$
\omega^2=1+u_T^2+u_\gamma^2,\quad u_\gamma=\phi\omega
$$
There are two cases about $|u_T|$. The first case is  $|u_T|^2(x_0)\leq 1$. From $(1-\phi^2)\omega^2=1+|u_T|^2$ we see that
$$
\omega(x_0)\leq \sqrt{\F{2}{1-\phi_0^2}}
$$
Thus the bound of $\omega$ is established and we are done. Otherwise we have the second case $|u_T|^2(x_0)\geq 1$. By the maximality of $\omega$, we observe that
$$
\omega_\gamma=u_\gamma u_{\gamma\gamma}+u_T u_{T\gamma}\leq 0
$$
and $\omega_T=0$. Consequently \eqref{ded} is simplified as
\begin{gather}
TT(u)=-\F{\phi\phi_T\omega^2}{u_T}\quad \gamma T(u)=\phi_T\omega+k u_T \label{eq:gammtu}\\
T\gamma(u)=\phi_T\omega
\end{gather}
By conclusion (ii) in Lemma \ref{lm:estd}, $\nb_\gamma \gamma=\nb_\gamma T=0$. Then  $u_{\gamma\gamma}=\gamma(\gamma(u)), u_{TT}=T(T(u))$. This gives us
\be\label{eq:gg}
\gamma(\gamma(u))u_\gamma + \gamma(T(u))u_T\leq 0
\ene
The main idea is to turn the above expression into one which contains only the first derivative of $u$. The key step is to solve $\gamma(\gamma(u))$. This is accomplished by the expedient of rewriting $g^{ij}u_{ij}$. By definition we have
\be\label{eq:com:ut}
\begin{split}
	&g^{ij}u_{ij} =g^{TT}(T(T(u))+g^{T\gamma}(T(\gamma(u))+g^{\gamma T}(\gamma (T(u))\\
	&+g^{\gamma\gamma}\gamma(\gamma(u))-g^{TT}\la \nb_TT, Du\ra-g^{T\gamma}\la \nb_T\gamma, Du\ra\\
	&-g^{\gamma T}\la \nb_{\gamma}T,Du\ra-g^{\gamma\gamma}\la \nb_{\gamma}\gamma, Du\ra
\end{split}
\ene
By conclusion (ii) in Lemma \ref{lm:estd}, the second line above becomes
$$
g^{\gamma\gamma}\gamma(\gamma(u))-k\F{(1+u_\gamma^2) u_\gamma}{\omega^2}-k\F{u_T^2u_\gamma }{\omega^2}=
g^{\gamma\gamma}\gamma(\gamma(u))-ku_\gamma
$$
With $u_\gamma =\phi\omega$ the first line turns into
\begin{align*}
&-(1+u_\gamma^2)\F{\phi\phi_T}{u_T}-2\F{\phi_T}{\omega} u_\gamma u_T-k\F{u_T^2u_\gamma}{\omega^2}\\
&=-(1+u_\gamma^2+2u_T^2)\F{\phi\phi_T}{u_T}-k\F{u_T^2u_\gamma}{\omega^2}
\end{align*}
Combining these expressions together, equation \eqref{eq:com:ut} is simplified as
$$
g^{ij}u_{ij}+(\omega^2+u_T^2)\F{\phi\phi_T}{u_T}+k u_\gamma \F{\omega^2+u_T^2}{\omega^2}-(1-\phi^2)\gamma\gamma(u)=0
$$
Multiplying the above equation by $u_\gamma$ and inequality \eqref{eq:gg} by $1-\phi^2$, summing them together and substituting $\gamma T(u)$ with equation \eqref{eq:gammtu} we obtain
\be\label{dded}
\begin{split}
	0 &\geq g^{ij}u_{ij}u_{\gamma}+(\omega^2+u_T^2)\F{\phi\phi_Tu_\gamma}{u_T}+(1-\phi^2)u_T\phi_T\omega\\
	&+k \{u^2_\gamma \F{\omega^2+u_T^2}{\omega^2}+u_T^2(1-\phi^2)\}
\end{split}
\ene
Next substituting $u_{\gamma}$ with $\phi\omega$ we observe the following identity:
\begin{gather*}
(\omega^2+u_T^2)\F{\phi\phi_Tu_\gamma}{u_T}+(1-\phi^2)u_T\phi_T\omega=\F{\phi_T\omega}{u_T}(\omega^2-1)\\
k \{u^2_\gamma \F{\omega^2+u_T^2}{\omega^2}+u_T^2(1-\phi^2)\}=k(\omega^2-1);
\end{gather*}
With them we continue to simplify inequality \eqref{dded} and obtain
\be\label{eq:sst}
\phi\omega g^{ij}u_{ij}+\F{\phi_T\omega}{u_T}(\omega^2-1)+k(\omega^2-1)\leq 0
\ene
Note that
$$
\omega^2-1=\F{\phi^2}{1-\phi^2}+\F{u^2_T}{1-\phi^2}
$$
This implies that $\F{|u_T|}{\sqrt{1-\phi^2}}\leq \omega$.
Therefore dividing inequality \eqref{eq:sst} by $\omega$ yields that
\be
k\omega- \omega\F{|\phi_T|}{\sqrt{1-\phi^2}}\leq \phi | g^{ij}u_{ij}| +\F{k}{\omega}+\F{\phi^2|\phi_T|}{(1-\phi^2)|u_T|}
\ene
By our assumptions and $|u_T|\geq 1$ at $x_0$,
$$
\delta_0\omega \leq ( k-\F{|\phi_T|}{\sqrt{1-\phi^2}})\omega \leq (C_0+k_0+\F{\phi_0k_0}{1-\phi_0^2})
$$
In either case we establish the bound of $\omega$. The proof is complete.
\ep
Next we prove the universal bound of $\omega$ under the setting of Theorem \ref{thm:MT3}. So we can complete conclusion (i) of Theorem \ref{thm:MT3}. In the following derivation the role of the assumption $Ric\geq 0$ is to guarantee that the maximum of $\omega$ has to occur on the boundary.
\bl\label{akd} Under the assumption in Theorem \ref{thm:MT3}, we have
$\omega\leq C_1$ where $C_1$ is a constant depending on $ k_0, \delta_0, \phi_0$ and the initial condition.
\el
\bp Using Theorem \ref{thm:bound_ut} in the case of $\psi(x,u, Du)\equiv 0$ and setting $\chi_0\equiv 0$, we
get
\be\label{guar}
|u_t|\leq \max_{x\in \bar{\Omega}}|u_t(x,0)|:=C_0
\ene
where $C_0$ is a constant independent of time.
Applying equation \eqref{ev:omega} into the case $\psi(x,u,Du)\equiv 0$ and the assumption $Ric\geq 0$ on $M^2$, we obtain
\begin{align}
L\omega &= (|A|^2+Ric(v_N, v_N))\omega+\F{2}{\omega}g^{il}\omega_i\omega_l \notag\\
&\geq |A|^2\omega+\F{2}{\omega}g^{il}\omega_i\omega_l \label{base:mp}\
\end{align}
where $L\omega=g^{ij}\omega_{ij}-\omega_t$.\\
\indent Fix any $t^*>0$. By inequality \eqref{base:mp}, the maximal principle of parabolic equations implies that the maximum of $\omega$ on $\bar{\Omega}\times [0, t^*]$  occurs at $(x_0, t_0)$ for $t_0=0$ or $x_0\in\P\Omega$. If in the first case, we are already done. Suppose $x_0\in\P\Omega$. Estimate \eqref{guar} can guarantee that on $\P\Omega$,
$$
| g^{ij}u_{ij}|=|u_t|\leq C_0
$$
For fixed $t_0$, this allows us to apply Lemma \ref{lm:convex:bound} and conclude
$$
\omega(x_0,t_0)\leq C_1=C_1(C_0,\delta_0, k_0, \phi_0);
$$
Thus for any $t\in [0,t^*]$, $\omega\leq C_1$ where $C_1$ is also a constant independent of $t$. Since $C_0$ relies on the initial condition, the proof is accomplished.
\ep
As in \cite{AW94}, we are interested in asymptotic behaviors of the following elliptic problem.
\be\label{eq:2:C}
\left\{\begin{split}
	&g^{ij}(Du)u_{ij}=C&\quad &\text{on $\Omega$}\\
	&u_\gamma =\phi(x)\omega&\quad & x\in\P\Omega
\end{split}\right.
\ene
where $C$ is a uniquely determined constant given by
\be\label{eq:dc}
C=-\F{\int_{\P\Omega}\phi(x)dx}{\int_{\Omega}\omega^{-1}dx}
\ene
A classical way to seek the solution to problem \eqref{eq:2:C} is to solve the following capillary problem:
\be\label{eq:2:CSc}
\left\{\begin{aligned}
	&g^{ij}(Du)u_{ij}=\Sc u\quad \text{on $\Omega$}\\
	\quad &u_\gamma =\phi(x)\omega\quad x\in \P\Omega
\end{aligned}\right.
\ene
for $\Sc>0$.
\bl\label{lm:ts} Under the assumptions of Theorem \ref{thm:MT3}, problem \eqref{eq:2:C} has a unique smooth solution $u(x)\in C^{\infty}(\bar{\Omega})$ (up to a constant).
\el
\br The existence of the solution to problem \eqref{eq:2:C} is highly nontrivial. We need four conditions: the positivity of Ricci curvature, the solvability of Capillary problems with positive gravity, the convex boundary of $\Omega$ and the dimension of $M^2$ is two.
\er
\bp Without loss of generality, we suppose $\Sc\in (0,1]$. By Theorem \ref{thm:MT22}, there is a $u_{\Sc}(x)\in C^{\infty}(\bar{\Omega})$ as the unique solution to problem \eqref{eq:2:CSc}. Corollary \ref{cor:est:u} implies that
\be\label{bound:Sc}
|\Sc u_{\Sc}|\leq c_0
\ene
where $c_0$ is a constant independent of $\Sc$. This implies $|g^{ij}(Du_\Sc)u_{ij}|\leq c_0$. According to equation \eqref{ev:omega}, $\omega$ satisfies that
\begin{align*}
g^{ij}\omega_{ij} &=(|A|^2+Ric(v_M, v_M))\omega+\F{2}{\omega}g^{il}\omega_i\omega_l\\
&\geq \F{2}{\omega}g^{il}\omega_i\omega_l
\end{align*}
where $Ric\geq 0$ from the assumptions of Theorem \ref{thm:MT3}. From the maximum principle the maximum of $\omega$ has to occur at the boundary. Again by the assumptions of Theorem \ref{thm:MT3} and \eqref{bound:Sc}, Lemma \ref{lm:convex:bound} implies that
$$
|Du_\Sc|\leq \omega \leq c_1(c_0, k_0, \delta_0,\phi_0)
$$
which also does not rely on $\Sc$. Therefore $D(\Sc u_{\Sc})$ converges to $0$ as $\Sc\rightarrow 0$. By inequality  \eqref{bound:Sc}, one sees that $\Sc u_{\Sc}$ converges to a constant as $\Sc\rightarrow 0$ (after choosing a suitable subsequence). Fix a point $x_0\in \Omega$, then $u_{\Sc}(x)-u_{\Sc}(x_0)$ is uniformly bounded since $|Du_\Sc|\leq c_1$. Then $u_{\Sc}(x)-u_{\Sc}(x_0)$ converges uniformly to a smooth solution to problem \eqref{eq:2:C}. We establish the existence part. \\
\indent Now we assume $f, g$ are two solutions to problem \eqref{eq:2:C} with some constant $C$. Adding some positive constant we can suppose $u=f-g\geq 0$. According to Lemma \ref{difference} in appendix C, $u$ satisfies $g^{ij}(Du_1)u+\tilde{b}_i u_i=0$. From the maximum principle, $u$ can achieve its nonnegative maximum at some $x_0\in \P\Omega$. By Lemma \ref{gamma:vanish} in appendix C, $u_\gamma(x_0)=0$. The Hopf lemma implies that $u$ is a constant. Thus the solution to problem \eqref{eq:2:C} is unique up to a constant.
\ep
We present two results  in \cite{AW94}. Their proofs only rely on the maximum principle of parabolic equations (see \cite{Lie96}). So the conclusions are still true for general Riemannian surfaces. \\
\indent  A direct observation is that for a solution $u_\infty(x)$ to problem \eqref{eq:2:C}, $u_\infty(x)+Ct$ solves problem \eqref{eq:2:ts}. Namely, $u_\infty(x)+Ct$ translates upwards (downwards) with speed $C>0(C<0)$. This gives an oscillation bound of solutions to problem \eqref{eq:2:ts}.
\bl \label{bound:ut}
Under the assumption of Theorem \ref{thm:MT3}, a solution $u(x,t)$ to problem \eqref{eq:2:ts} satisfies $$
|u(x,t)-Ct|\leq c_2
$$
for some constant $c_2$ depending on $\Omega$. Here $C$ is unique determined by $\Omega$ in Lemma \ref{lm:ts}.
\el
\bp The proof is exactly the same as Corollary 2.7 in \cite{AW94}.
\ep
Consequently all solutions to problem \eqref{eq:2:ts} on a fixed convex domain ( possible with different initial conditions) have the same asymptotic behavior with fixed prescribed contact angle (the same $\phi(x)$).
\bl\label{eqst} Suppose $\Omega$ satisfies the assumptions of Theorem \ref{thm:MT3}. Let $u_1$ and $u_2$ be two solutions to problem \eqref{eq:2:ts} with prescribed contact angle. Define $u:=u_1-u_2$. Then $u$ becomes a constant function as $t\rightarrow \infty$. In particular $u_\infty(x)$ is the solution of problem \eqref{eq:2:C} with a constant $C$, then $u_1(x,t)-u_\infty(x)-Ct$ converges uniformly to a constant.
\el
\bp The proof is exactly the same as Theorem 3.1 in \cite{AW94}.
\ep

Now we are ready to conclude Theorem \ref{thm:MT3}.
\bp(The proof of Theorem \ref{thm:MT3}) Conclusion (i) follows from Lemma \ref{akd}. Conclusion (ii) follows from Lemma \ref{eqst}. The proof is complete.
\ep
\section{\textit{Acknowledgements}}
The author is grateful to Prof. Zheng Huang for his discussions and encouragements. He also expresses his gratitude to the support of Prof. Lixin Liu.
He would also like to express sincere gratitudes to the referee for many constructive suggestions and helpful comments.
\appendix
\section{A counterexample}
\begin{figure}[h]
	\centering                       
	  \includegraphics[width=.5\textwidth]{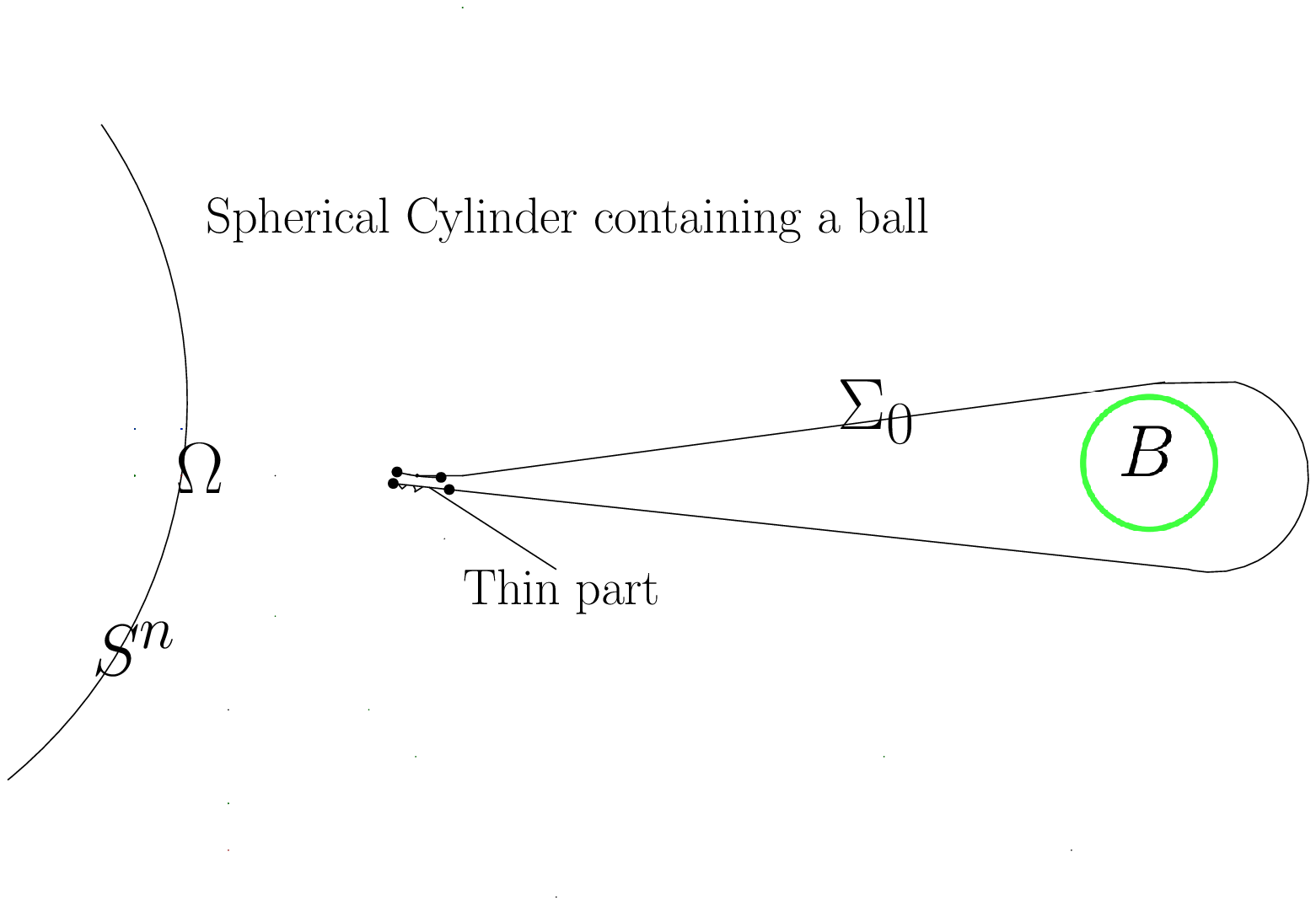}
	\caption{Spherical Cylinder}\label{fig:cylinder}
\end{figure}

\begin{figure}[ht]
	\centering                    
	\includegraphics[width=.5\textwidth]{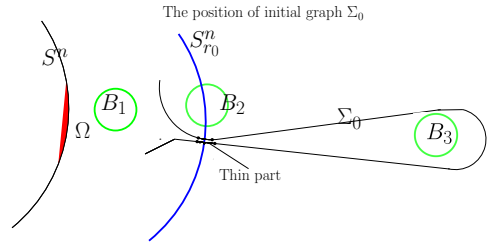}                         %
	\caption{Initial graph $\Sigma_0$}  \label{fig:two}                                  %
\end{figure}

This section is devoted to the proof of Theorem \ref{thm:example}. We will construct a graph in a warped product manifold such that its mean curvature flow only exists for finite time and loses the graphic property.  The content in this section is independent of other parts in this paper. Our purpose is to illustrate that Theorem \ref{thm:MT1} is optimal.  \\
\indent The warped product manifold we will discuss is given as follows.  Fix $n\geq 2$. Let $\R_b$ denote a complete warped product manifold given by
\be\label{ee:wp}
\R_b=(S^n\times (1,\infty), dr^2+\rho^2(r)\sigma_n)
\ene
where $(S^n,\sigma_n)$ is the unit sphere in $\R^{n+1}$ and $\rho(r):\R\rightarrow (0,\infty)$ is a smooth positive function such that $\rho'(0)=0$, $\rho''(r)\geq 0$ for all $r$ and $\rho(r)=r$ for $r\geq 1$. In this setting the set $\{(x,r):r\geq 1\}$ in $\R_b$ is a Euclidean space removed the unit ball. A mean curvature flow in this set is indeed a mean curvature flow in Euclidean space.\\
\indent A hypersurface in $\R_b$ is called as a \emph{(starshaped) graph} over a domain $\Omega$ in $S^n$ if it is a set $\{(x, u(x))\in \R_b: x\in \Omega \}$ where $u(x)$ is a smooth function on $\Omega$. \\
\indent The goal of this section is to show
\bt\label{thm:example} Given any bounded domain $\Omega$ in $S^n$, we can find a function $u_0(x)$, $u_0(x)>>1$, on $\Omega$ such that the mean curvature flow of its graph $\Sigma_0$ in \eqref{ee:wp} will lose its graphical property before time $t=1$ and form singularity.
\et
Before the proof, we need three preliminary facts. All of them are valid in the Euclidean part of $\R_b$. The first one is the well-known clearing-out lemma (Lemma 6.3 in \cite{Bra78}).
\bt \label{thm:cout} Let $n\geq 2$ and $\mathcal{H}^{n}$ be the $n$ dimensional Hausdorff measure in $\R^{n+1}$. Consider a smooth mean curvature flow of hypersurfaces $\{F_{t}(\Sigma)\}_{t\geq0}$ in $\R^{n+1}$. There exists a constant $c>0$ such that if
$\mathcal{H}^n(\Sigma\cap B^n_{\rho}(x_0))\leq \Sc\rho^2$ for some $x_0\in \R^{n+1}$, $\rho >0$ and $\Sc >0$, then
$\mathcal{H}^n(F_{t}(\Sigma)\cap B^n_{\F{\rho}{4}}(x_0))=0$ for $t=c\Sc\rho^2$. Here $B^n_r(x)$ is the Euclidean ball with radius $r$ centered at $x$.
\et
More discussion in detail can be found in \cite{Eck91}, appendix E in \cite{Eck04} and appendix E in \cite{Unt99}.\\
\indent The second one is the disjoint principle by \cite{Hui86}.
\bt Let $\Sigma^1_t$ be a mean curvature flow of a closed hypersurface and $\Sigma^2_t$ be a mean curvature flow of a compact hypersurface with some boundary condition. If they are disjoint at $t=0$, then $\Sigma^1_t$ are disjoint with the interior of $\Sigma^2_t$ for any time that both of them exist smoothly.
\et
\bp The proof is the same as Lemma 3.2 in \cite{Hui86}.
\ep
\indent  Let $\Omega$ be a domain in $S^n$. $\Omega\times (a,b)$ is called as a \textit{\textbf{spherical cylinder}} in $\R_b$ if $a > 1$. In terms of polar coordinates it is a set in $\R^{n+1}$ consisting of all points $X$ such that $|X|\in (a,b)$ and $|X|^{-1}X\in \Omega$.\\
\indent The third fact we need is stated as follows.
\bt\label{thm:sc} Given any domain $\Omega$ in $S^n$ and any positive number $r$, there is a $k_0>0$ such that for all $k\geq k_0$ the spherical cylinder $\Omega\times(k, 2k)$ contains a Euclidean ball with radius $r$ in $\R_b$. \et
Here the size of $\Omega\subset S^n$ may be very tiny. The rough idea is depicted in Figure \ref{fig:cylinder}.
\bp Now assume all $r$-coordinates are greater than one. We always work in the Euclidean part of $\R_b$.\\
\indent Let $S\Omega_k$ be the spherical cylinder $\Omega\times(k, 2k)$ where $k$ is determined later. It suffices to show that there is a point $p$ in $S\Omega_k$ such that its Euclidean distance to the boundary is greater than $r$. \\
\indent Let $d_{S^n}$ be the distance on the unit sphere. Fix a point $x_0\in \Omega$. Then there is a $\Sc_0$ such that $d_{S^n}(x_0,\P\Omega)\geq \Sc_0$. We denote the point $(x_0,\F{3k}{2})$ in $\R_b$ by $p$. Our goal is to show that this point gives the desirable property if $k$ is sufficiently large. \\
\indent In the following we compute the distance from $p$ to the boundary of $S\Omega_k$. The boundary of $S\Omega_k$, $\P S\Omega_k$, is composed of $\Sigma_1$ and $\Sigma_2$ given as follows:
\be
\Sigma_1=\{(x,r)\in \R_b: x\in \Omega, r=k, 2k\}\quad \Sigma_2=\{(x,r)\in \R_b: x\in \P\Omega, k\leq r\leq 2k\}
\ene
Let $d$ be the distance of $\R_b$. The definition of the warped product metric in \eqref{ee:wp} yields that
\begin{gather}
d(p, \Sigma_1)\geq \inf_{(x,r)\in \Sigma_1}|r-\F{3k}{2}|\geq \F{k}{2}\\
d(p, \Sigma_2)\geq k^2\inf_{x\in \P\Omega} d_{S^n}(x,x_0)=k^2\Sc_0.
\end{gather}
Therefore the distance of $p$ to the boundary of $S\Omega_k$ is 
$$
d(p,\P S\Omega_k)= k\min\{\F{1}{2},\Sc_0 k\}
$$
Letting $k$ go to $\infty$ we have $d(x,\P S\Omega_k)\rightarrow \infty$. This gives a $k_0$ such that for all $k\geq k_0$, $ d(p,\P S\Omega_k)\geq r$. One can easily observe that $d$ is just the Euclidean distance if $k_0$ is large enough. We complete the proof.
\ep  
Now we begin to the construction of $u_0(x)$ in Theorem \ref{thm:example}.
\bp(The proof of Theorem \ref{thm:example}) In the following all $r$-coordinates are assumed to be greater than one. We denote the graph of $u_0(x)$ by $\Sigma_0$. We fix two points $x_1, x_2$ in $\Omega$. Choose $d_0>0$ such that $d_0\leq \F{d_{S^n}(x_1, x_2)}{4}$ and
\be \label{def:omega1}
\Omega_1=\{x\in S^n: d_{S^n}(x,x_1)\leq d_0\}
\ene
is contained in $\Omega$.\\
\indent Let $c, \Sc$ be the constants in Theorem \ref{thm:cout}. By Theorem \ref{thm:sc}, there exists a sufficiently large number
\be\label{chose:k}
k_0\geq \max\{5\sqrt{2n}, \F{5}{\sqrt{c\Sc}}\}
\ene such that there are two Euclidean balls $B_1$ and $B_2$ with radius $\sqrt{2n}$ belonging to sphere cylinders $\Omega_1\times (k_0, 2 k_0)$ and $\Omega_1\times (3k_0, 4k_0)$ respectively. In $\R_b$ ( Euclidean space) the mean curvature flow of a sphere with Euclidean radius $\sqrt{2n}$ is a family of spheres with Euclidean radius $\sqrt{2n-2nt}$ and exists smoothly until at time $t=1$ it contracts to its center.\\
\indent Suppose $r_0$ is the $r$-coordinate of the center of $B_2$. Then $ 4k_0 \geq r_0\geq 3k_0$. Let $S^n_{r_0}$ be the set $\{(x,r_0):x\in S^n\}$ in $\R_b$. It is indeed the sphere with radius $r_0$ in Euclidean space $\R^{n+1}$. To approach our objective, we need $u_0(x)$ on $\Omega$ to fulfill at least the following two conditions.
\begin{itemize}
	\item[(p1)] The graph of $u_0(x)$ in $\R_b$ is disjoint with $B_1, B_2$.
	\item[(p2)] Let $\mathfrak{S}$ denote the intersection of $\Sigma_0$ and the spherical cylinder $S^n\times (r_0- \F{1}{\sqrt{c\Sc}}, r_0+\F{1}{\sqrt{c\Sc}})$ in $\R_b$. The Hausdorff measure of $\mathfrak{S}$ is smaller than $\F{1}{c}$ where $c$ is the constant from Theorem \ref{thm:cout}.
\end{itemize}
To this end, we consider a strictly decreasing positive function $v(r)$ for $r\geq 0$ as follows.
\be\label{def:v}
v(r)=\left\{\begin{aligned}
	&2k_0 \quad r\geq 2d_1\\
	& 6k_0  \quad r= d_1\\
	&k_1\quad r= \F{d_1}{2}\\
	&2k_1\quad r\leq \F{d_1}{4}
\end{aligned}\right.
\ene
with $|v'(r)|\leq \F{4k_0}{d_1}$ for $r\in (d_1, 2d_1)$. Here $k_1$ and $d_1$ will be determined lately where $d_1$ is sufficiently small.
Now we define a smooth function $u_0(x)$ as follows
\be
u_0(x)=v(d_{S^n}(x,x_2))
\ene
Then we assume that $d_1$ is sufficiently small such that $d(x,x_2)\geq 2d_1$ for all $x\in \Omega_1$. Thus  $u_0(x)\equiv 2k_0$ for $x\in \Omega_1$. As a result the graph $\Sigma_0$ satisfies property (p1).\\
\indent The property (p2) is realized to continue to choose sufficiently small $d_1$ as follows. According to assumption \eqref{chose:k} and $r_0\in [3k_0,4k_0]$, we observe that
$$
2k_0 <r_0-\F{1}{\sqrt{c\Sc}}<r_0+\F{1}{\sqrt{c\Sc}}< 6k_0
$$
Examining definition \eqref{def:v} we see that given any point $p\in \mathfrak{S}$ the $x$-coordinate of $p$ is contained in a set $E$ defined by $$
\{x\in \Omega: d_1\leq d_{S^n}(x, x_2)\leq 2d_1\}$$
Moreover we have
$$
2k_0\leq   |u_0(x)|\leq 6k_0,\quad  |Du_0(x)|\leq \F{4k_0}{d_1}
$$
when $x\in E$.
By Lemma \ref{area} in appendix B, the area of $\mathfrak{S}$ over $E$ is bounded above by
$$
\int_{E}\sqrt{1+\F{|Du_0(x)|^2}{u_0(x)^2}} u^n_0(x)dx\leq (6k_0)^n\sqrt{1+\F{4}{d_1^2}} \int_{E}dx =O(d_1^n)\sqrt{1+\F{4}{d_1^2}}
$$
Here we use the fact the dimension of $S^n$ is $n$. Since $n\geq 2$, $O(d_1^n)\sqrt{1+\F{4}{d_1^2}}$ converges to $0$ as $d_1\rightarrow 0$. As a result we can choose a sufficiently small $d_1$ such that the area of $\mathfrak{S}$ is less than $\F{1}{c}$. This is the property (p2).\\
\indent  Now only $k_1$ is left to be determined. We define a domain $\Omega_2$
$$
\Omega_2=\{x\in\Omega: d_{S^n}(x_2, x)\leq d_2\}
$$
where $d_2\leq \min\{d_0, \F{d_1}{4} \}$ such that $\Omega_2$ is contained in $\Omega$ and is disjoint with $\Omega_1$.
Applying Theorem \ref{thm:sc}, there exists a $k_1>6k_0$ such that the spherical cylinder $\Omega_2\times (k_1, 2k_1)$ contains a Euclidean ball $B_3$ with radius $\sqrt{2n}$. Fix this $k_1$. From the definition in \eqref{def:v}, we conclude
\begin{itemize}
	\item[(p3)]  the graph of $u_0(x)$ is disjoint with the ball $B_3$ with radius $\sqrt{2n}$. The $r$-coordinate of the center of $B_3$ is less than $2k_1$ and greater than $k_1$.
\end{itemize}
Finally we complete the construction of $u_0(x)$ with the property (p1), (p2) and (p3). Moreover
$u_0(x)\equiv 2k_0$ on $\Omega_1$ and $u_0(x)\equiv 2k_1$ on $\Omega_2$. The rough shape of $\Sigma_0$ is described in Figure \ref{fig:two}.
\\
\indent Suppose the mean curvature flow of $\Sigma_0$, noted as $\Sigma_t$, exists smoothly on $[0,1]$ and preserves the graphical property.
\begin{itemize}
	\item [(i)]At time $t=0$, $u_0(x)\equiv 2k_0$ on $\Omega_1$. By property (p1) the disjoint principle implies there is at least one point $p_1$ on $\Sigma_1$ such that its $r$-coordinate is strictly less than $r_0$ at time $t=1$.
	\item[(ii)] Similarly, by property (p2), there is at least one point $p_2$ on $\Sigma_1$ such that its $r$-coordinate is greater than $k_1$ at $t=1$.
\end{itemize}
Note that $r_0\in (3k_0, 4k_0)$ and $k_1\geq 6k_0$. By the continuity, there is one point $p_3$ on $\Sigma_1$ such that its $r$-coordinate is $r_0$. At time $t=0$, property (p2) implies that
$$
\mathcal{H}^n(\Sigma_0\cap B^n_{\F{1}{\sqrt{c\Sc}}}(p_3)))\leq \mathcal{H}^n(\mathfrak{S})\leq \F{1}{c}
$$
Applying Theorem \ref{thm:cout}, $p_3$ does not belong to $\Sigma_1$. This is the only but essential place we need the clearing out lemma Theorem \ref{thm:cout}. In fact the same conclusion is also valid for any point on $S_{r_0}^n$. Hence this gives a contradiction. \\
\indent This implies that the mean curvature flow of $\Sigma_0$ forms its singularity before time $t=1$. Moreover it loses the graphic property. Otherwise, the $r$-coordinate of the flow satisfies a parabolic equation generated by
$
(\F{d F}{dt})^{\bot}=\vec{H}
$
where $F(x,t)=(x,u(x,t))$ and $\vec{H}$ is the mean curvature vector of $F(x,t)$. It is impossible for the singularity to occur.  We skip the derivation about the expression of the parabolic equation. It is just a direct and tedious computation. The proof is complete.
\ep

\section{An area formula in warped product manifolds}
We give a formula for computing the area of starshaped graphs in warped product manifolds.  A warped product manifold is defined by
\be
(M\times I, dr^2 +\rho(r)^2 \sigma)
\ene
where $(M,\sigma)$ is a Riemannian manifold, $I$ is an open interval of $\R$ and $\rho$ is a smooth positive function on $I$. Let $\Omega$ denote a bounded domain and $\Sigma$ be the set $ \{(x,u(x)\}$ where $u(x)$ is a smooth function on $\Omega$.
\bl\label{area}
The area of $\Sigma$ is given by
\be
\int_{\Omega} \sqrt{1+\F{|Du|^2}{\rho^2(u(x))}}\rho^{n}(u(x))dvol
\ene
where $dvol$ is the volume form of $(M,\sigma)$.
\el
\bp Fix a point $p=(x,u(x))$. We choose a local frame $\{\P_i\}$ near $x$ with the property that $\la \P_i, \P_j\ra=\delta_{ij}$ and $u_i=0$ for $i\geq 2$ at $x$. We obtain a local frame on $\Sigma$ near $p$ given as
$$
\{ \P_i+u_i\P_r \}_{i=1}^n
$$
Now all expressions are evaluated at $p=(x,u(x))$. Under our assumptions
$$
|Du(x)|^2=|u_1|^2, dvol=dx^1\wedge \cdots\wedge dx^n
$$
Thus the volume form on $\Sigma$ takes the form $\sqrt{det(\la X_i, X_j\ra)}dx^1\wedge \cdots \wedge dx^n$
where $X_i$ denotes $\P_i+u_i\P_r$.  The conclusion follows from
$$
det(\la X_i, X_j\ra)=\rho^{2n}(u(x))(1+\F{|Du|^2}{\rho^2(u(x))})
$$
since $\la X_i, X_j \ra=\rho^2(u(x))\delta_{ij}+u_iu_j$ and $u_i=0$ for $i\geq 2$.
\ep
\section{Some results in Partial Differential Equations}
In this section, we recall some technique results in partial differential equations. We define an operator,
$$
Lu= a^{ij}(Du)u_{ij}+b_i(Du)u_i-\phi(x,u, Du)-k\P_t u
$$
where $(a^{ij})$ is a smooth, positive definite matrix and $k$ is a constant equal to $1$ or $0$.
\bl\label{difference} Suppose $u_i(x)$ are smooth solutions of $Lu=C_i$ on a bounded smooth domain $\Omega$ for $i=1,2$. Then $u:=u_1-u_2$ satisfies
$$
a^{ij}(u_1)u_{ij}+\tilde{b}_i u_i-\phi_u(x, *, Du)u-k\P_t u=C_1-C_2
$$
for some smooth function $\tilde{b}_i$. Here $*$ denotes some unknown function of $x$.
\el
\bp First we have $Lu_1-Lu_2\equiv C_1-C_2$. One easily observes that
\begin{align*}
&Lu_1-Lu_2= a^{ij}(Du_1)(u_1)_{ij}-a^{ij}(Du_2)(u_2)_{ij}\label{c:B}\\
&+b_i(Du_1)(u_1)_i-b_i(Du_2)(u_2)_i\\
&-\phi(x,u_1, Du_1)+\phi(x,u_2, Du_2)-k\P_t u
\end{align*}
We mainly apply the Taylor formula.
The right hand of the first line may be written as
\begin{align*}
&a^{ij}(Du_1)(u)_{ij}+ (a^{ij}(Du_1)-a^{ij}(Du_2))(u_2)_{ij}\\
&=a^{ij}(Du_1)(u)_{ij}+a^{ij}_k(*)(u_2)_{ij}u_k
\end{align*}
We compute the second line with a similar idea and obtain that
$$
b_i(Du_1)(u_1)_i-b_i(Du_2)(u_2)_i=(b_i(*))_k( u_1)_iu_k+ b_k(Du_2)u_k
$$
As for the last line, we see that
\begin{align*}
&-\phi(x,u_1, Du_1)+\phi(x,u_2, Du_2)-k\P_t u\\
&=-\phi_u(x,*, Du_1)u+\phi_{u_k}(x,u_2, *)u_k-k\P_t u
\end{align*}
Putting these facts together and combining like terms we obtain the equation at the beginning.
\ep
\bl\label{gamma:vanish} Assume $f, g$ are two smooth functions on $\bar{\Omega}$ satisfying
\be\label{ab}
f_\gamma=\phi(x)\sqrt{1+|Df|^2}\quad   g_\gamma=\phi(x)\sqrt{1+|Dg|^2}
\ene
on $\P\Omega$. If $u:=f-g$ attains its maximum or minimum at $x_0\in \P\Omega$, then $u_\gamma(x_0)=0$.
\el
\bp Assume $u:=f-g$ attains its maximum or minimum at $x_0\in \P\Omega$, then $f_T(x_0)=g_T(x_0)$ for any tangent vector $T$ on $\P\Omega$. Let $\{T_i, \gamma\}$ be an orthonormal unit frame at $x_0$.  Thus at $x_0$ we have
$$
|Df|^2 =|f_{T_i}|^2+|f_\gamma|^2,  |Dg|^2 =|g_{T_i}|^2+|g_\gamma|^2,
$$
We denote $|f_{T_i}|^2=|g_{T_i}|^2$ by $a$. By equation \eqref{ab} we observe that
$$
\F{f_\gamma}{\sqrt{1+a^2+f_\gamma^2}}= \F{g_\gamma}{\sqrt{1+a^2+g_\gamma^2}}
$$
The monotonicity of $\F{q}{\sqrt{1+a^2+q^2}}$ as a function of $q$ yields $f_\gamma=g_\gamma$. The proof is complete.
\ep
\section{The proof of Lemma \ref{lm:cn}}
Lemma \ref{lm:cn} shall be already known in \cite{GB96},\cite{Hui89A} and \cite{AW94}. It is a sufficient condition on the uniformly convergence of the solution in problem \eqref{eq:first}. For convenience of the reader, we include its proof here.
\bp Suppose $u_t(x,t)$ does not converges uniformly to 0. There are a constant $\Sc_0>0$, a sequence $\{x_n\}\in \bar{\Omega}$ and $\{t_{n}\}\rightarrow \infty$ such that
$$
|u_t(x_n, t_n)|\geq \Sc_0
$$ Without loss of generality, we assume $t_{n+1}-t_{n}\geq 1$. Note that $\bar{\Omega}$ is bounded. By assumptions and the classical Schauder estimates of elliptic equations, all derivatives of $u(x,t)$ are uniformly bounded. Thus after suitable choosing subsequence, there is a point $x_0\in \bar{\Omega}$ and $\Sc<\F{1}{2}$ such that
$$ |u_t(x_0,t)|\geq \F{\Sc_0}{2} $$
for all $t\in (t_n-\Sc, t_n+\Sc)$ where $n$ is sufficiently large. This in turn implies that there is  a compact subdomain $\Omega_{x_0}$ containing $x_0$ such that for all $x\in \Omega_{x_0}$ and $t\in (t_n-\Sc, t_n+\Sc)$, we have
$$ |u_t(x,t)|\geq \F{\Sc_0}{4} $$
again because all higher derivatives of $u(x,t)$ are uniformly bounded.
Thus
$$
\int_{\Omega_{x_0}}\int_{t_n-\Sc}^{t_{n}+\Sc}|u_t|^2 dtdx\geq \F{\Sc_0}{2}vol(\Omega_{x_0})2\Sc>0
$$
as $t_n\rightarrow \infty$. This gives a contradiction to the fact $\int_0^T\int_{\omega}|u_t|dxdt$ is uniformly bounded for any $T$. Hence $u_t(x,t)$ converges uniformly to 0. \\
\indent Suppose $u(x,t)$ does not converges uniformly to a smooth function $u_\infty(x)$. There are a constant $\Sc_1>0$, a sequence $\{y_n\}\in \bar{\Omega}$ and $\{s_n\}\rightarrow \infty$ such that
\be\label{eq:sde}
|u(y_n, s_{2n-1})-u(y_n, s_{2n})|\geq \Sc_1
\ene
Since $\bar{\Omega}$ is compact and all derivatives of $u(x,t)$ are uniformly bounded, we can assume there is a sequence $\{n_k\}_{k=1}^\infty$ such that $u(x, s_{2n_k})$ and $u(x,s_{2n_{k}-1})$ converge to $u_{1,\infty}(x)$ and  $u_{2,\infty}(x)$ respectively. Since $u_t$ converges uniformly to $0$, $u_{i,\infty}(x)$ is the solution to problem \eqref{eq:second} with $\psi(x,u,Du)$
for $i=1,2$. Moreover $u_{1,\infty} \neq u_{2,\infty}$ by inequality \eqref{eq:sde}. \\
\indent We will see a contradiction case by case. If $\psi(x,u,Du)=h(x,u)Du$ with $h_u(x,u)$ greater than $h_0 >0$, then Lemma \ref{lm:unique} says that there is at most one solution for problem \eqref{eq:second} with $\psi(x,u,Du)$. Hence this is a contradiction since $u_{1,\infty} \neq u_{2,\infty}$. If $\psi(x,u, Du)\equiv 0$, problem \eqref{eq:second} only has constant solutions. Hence $u_{i,\infty}\equiv C_i$ for $i=1,2$ and we can assume $C_1>C_2$. There exists a $t_0$ such that  $\min_{x\in\bar{\Omega}}u(x, t_0)> C_2$. Applying the maximal principle to problem \eqref{eq:first} on $\Omega\times[t_0, \infty)$, we get for all $t>t_0$ $u(x,t)\geq \min_{x\in\bar{\Omega}}u(x, t_0)> C_2$. This is impossible because $C_2$ is the limit of $u(x, s_{2n_k-1})$ where $s_{2n_k -1}$ goes to infinity. Therefore when $\psi(x,u,Du)\equiv 0$, then $u(x,t)$ converges uniformly to a constant. The proof is complete.
\ep
\bibliographystyle{abbrv}	
\bibliography{Ref_Thesis}
\end{document}